\magnification 1200
\newread\rfaux 
\openin\rfaux=\jobname.AUX 
\ifeof\rfaux \message{Can't find \jobname.AUX!!!}
   \else \closein\rfaux 
     \input \jobname.AUX \fi 
 
\newwrite\wfaux 
\immediate\openout\wfaux=\jobname.AUX 
 
\newcount\eqnum  \newif\ifeqnerr 
 
\def\neweq#1{{\global\advance\eqnum by 1 
    \edef\chqtmp{(\the\eqnum)}
    \immediate\write\wfaux{\string\def \string#1{\chqtmp}}
    \checkeq#1}} 
 
\def\checkeq#1{\ifx#1\chqtmp\global\eqnerrfalse 
    \else \global\eqnerrtrue 
       \ifx#1\zundefined 
           \def\zz{NOT DEFINED}\global\let#1\chqtmp 
       \else \def\zz{preassigned #1}\fi 
        \message{ALLOCATION ERROR for \string#1: 
             \zz, should be \chqtmp.}\fi} 
 
\def\chq#1{\neweq#1\ifeqnerr***\fi#1} 
\def\cheqno{\eqno\chq} 
 
\def\pagemac#1{
    \write\wfaux{\string\def \string#1{\the\pageno}}
    \ifx#1\zundefined 
      \message{Page macro \string#1 NOT DEFINED.}\gdef#1{***}\fi} 
 
 
\def\neweqx#1#2{{\edef\chqtmp{(\the\eqnum#1)}
    \immediate\write\wfaux{\string\def \string#2{\chqtmp}}
    \checkeq#2}}

\font\tenmsb=msbm10   
\font\sevenmsb=msbm7
\font\fivemsb=msbm5
\newfam\msbfam
\textfont\msbfam=\tenmsb
\scriptfont\msbfam=\sevenmsb
\scriptscriptfont\msbfam=\fivemsb
\def\Bbb#1{\fam\msbfam\relax#1}
\let\nd\noindent 
\def\qed{\hbox{\hskip 6pt\vrule width6pt height7pt depth1pt \hskip1pt}}
\def\natural{{\rm I\kern-.18em N}}

\def\G{{\cal G}}
\def\B{{\cal B}}   
\def\E{{\Bbb E}}
\def\H{{\Bbb H}}
\def\F{{\cal E}}
\def\O{{\cal O}}
\def\P{{\cal P}}        
\def\S{{\Bbb X}}

\def \T{{\cal T}}  
\def \t{{\rm top}} 

\def\integer{{\rm Z\kern-.32em Z}}
\def\chix{{\raise.5ex\hbox{$\chi$}}}
\def\lowR{{\lower.5exR}}
\def\Z{{\Bbb Z}}
\def\real{{\rm I\kern-.2em R}}
\def\R{{\Bbb R}}
\def\Q{{\Bbb Q}}
\def\complex{\kern.1em{\raise.47ex\hbox{
            $\scriptscriptstyle |$}}\kern-.40em{\rm C}}

\def\vs#1 {\vskip#1truein}
\def\hs#1 {\hskip#1truein}

\def\Month{\ifcase\number\month \relax\or January \or February \or
  March \or April \or May \or June \or July \or August \or September
  \or October \or November \or December \else \relax\fi }
\def\date{\Month \the\day, \the\year}

  \hsize=6truein        \hoffset=.25truein 
  \vsize=8.8truein      
  \pageno=1     \baselineskip=12pt
  \parskip=0 pt         \parindent=20pt
  \overfullrule=0pt     \lineskip=0pt   \lineskiplimit=0pt
  \hbadness=10000 \vbadness=10000 
\pageno=0

\footline{\ifnum\pageno=0\hss\else\hss\tenrm\folio\hss\fi}
\hbox{}
\vskip 1truein\centerline{{\bf UNIQUENESS AND SYMMETRY}}
\vskip .1truein\centerline{{\bf IN PROBLEMS OF OPTIMALLY DENSE PACKINGS}}
\vskip .2truein\centerline{by}
\vskip .2truein
\centerline{Lewis Bowen,${}^1$ 
\footnote*{Research supported in part by NSF Vigre Grant DMS-0135345}
\ \ Charles Holton,${}^2$
\footnote{**}{Research supported in part by NSF Vigre Grant DMS-0091946}}
\vs.1
\centerline{{Charles Radin${}^2$ 
\footnote{***}{Research supported in part by NSF Grant DMS-0071643 and
\hfill\break \indent Texas ARP Grant 003658-158\hfil}}\ \ 
and Lorenzo Sadun${}^2$}
\vskip .2truein\centerline{\vbox{
${}^1$\ \ Department of Mathematics, University of California at Davis, CA
\vskip.1truein
${}^{2}$\ \ Department of Mathematics, University of Texas, Austin, TX}}

\vs.5
\centerline{{\bf Abstract}}
\vs.1 \nd
We analyze the general problem of determining optimally dense
packings, in a Euclidean or hyperbolic space, of congruent copies of
some fixed finite set of bodies. We are strongly guided by examples of
aperiodic tilings in Euclidean space and a detailed analysis of a new
family of examples in the hyperbolic plane. Our goal is to understand
qualitative features of such optimum density problems, in particular
the appropriate meaning of the uniqueness of solutions, and the role
of symmetry in classfying optimally dense packings.
\vs.8
\vs.2
\centerline{Subject Classification:\ \ 52A40, 52C26, 52C23}
\vfill\eject

\nd {\bf I. Introduction}
\vs.1
The objects of our study are the densest packings, particularly of
balls and polyhedra, in a space of infinite volume. These optimization
problems are of a fundamental nature and have been of interest for
many years; for a survey see the classic texts [Feje] and [Roge], and
the review [FeKu]. They appear as part of the eighteenth problem in
Hilbert's list [Hilb], a problem devoted to fundamental domains in
spaces of constant curvature, and the manner in which the domains fill
up such a space.

Most interest has centered on densest packings in the Euclidean spaces
$\E^n$, notably when the dimension $n$ is 2 or 3, but we will see that
analysis of the problem in hyperbolic spaces $\H^n$ can clarify some
issues for the Euclidean problems so we consider the more general
problem in the $n$ dimensional spaces $\S^n$, where $\S^n$ will stand
for either $\E^n$ or $\H^n$. (It would be reasonable to generalize our
considerations further, to symmetric spaces, and even to include
infinite graphs, but as we have no noteworthy results in that
generality we felt it would be misleading to couch our considerations
in that setting.)

Let us begin with some notation and basic features of density. Given
some finite collection $\B$ of ``bodies'' in $\S^n$ -- a body being a
connected compact set with dense interior and boundary of volume 0 
-- we consider ``packings''
of $\S^n$ by the bodies: that is, collections $P$ of congruent copies
of the bodies in which interiors of bodies do not intersect. Denoting
by $B_r(p)$ the closed ball in $\S^n$ of radius $r$ and center
$p$, we define the ``density relative to $B_r(p)$'' of a packing
$P$ as:

$$D_{B_r(p)}(P)\equiv {\sum_{\beta \in P} m_{\S^n}[\beta \cap B_r(p)]
\over m_{\S^n}[B_r(p)]}, \cheqno\eqdensrho $$

\nd where $m_{\S^n}$ is the usual measure on $\S^n$. Then, assuming the
limit exists, we define the ``density'' of $P$ as:

$$D(P)\equiv \lim_{r \to\infty}
{\sum_{\beta \in P} m_{\S^n}[\beta \cap B_r(p)]\over m_{\S^n}[B_r(p)]}. 
\cheqno\eqdens $$

It is not hard to construct packings $P$ for which the limiting
density $D(P)$ does not exist, for instance by the adroit choice of
arbitrarily large empty regions so that the relative density
oscillates with $r$ instead of having a limit. (In hyperbolic space
the limit could exist but depend on $p$, which we also consider
unacceptable.) The possible nonexistence of the limit of \eqdens\ is
an essential feature of analyzing density in spaces of infinite
volume; density is inherently a {\it global} quantity, and
fundamentally requires a formula somewhat like \eqdens\ for its
definition [Feje], [FeKu]. We discuss this further below.

Certainly the most important examples for which we have significant
information about the densest packings are the densest packings of
balls of fixed radius in $\E^n$ for $n=2$ and 3. (For a recent survey
of this problem in higher dimensions see [CoGS]). It will be useful in
discussing these problems to make use of the notion of ``Voronoi
cell'', defined for each body $\beta$ in a packing $P$ as the closure
of the set of those points in $\S^n$ closer to $\beta$ than to any
other body in $P$. A noteworthy feature of the $n=2$ example is, then,
that in the optimal packing (see Figure 1) the Voronoi cell of every
disk (the smallest regular hexagon that could contain the disk) has
the property that the fraction of the area of this cell taken up by
the disk is strictly larger than for any other Voronoi cell in any
packing by such disks. (Intuitively, the optimal configuration is
simultaneously optimal in all local regions.) As for $n=3$, it is
generally felt that the densest lattice packing (i.e., the face
centered cubic) achieves the optimum density among all possible
packings, along with all the other packings made by layering
hexagonally packed planar configurations, such as the hexagonal close
packed structure; see [Roge]. 
There are claims in the literature by Hsiang [Hsia]
and by Hales [Hale] for proofs of this, and there is hope that the
problem will soon be generally accepted as solved.

Less well known but perhaps next in significance as examples of
optimal density (see [Miln]) are the various ``aperiodic tilings'',
especially the ``Penrose kite \& dart tilings'', the tilings of $\E^2$
by congruent copies of the two polygons of Figure 2. (A portion of a
kite \& dart tiling is shown in Figure 3.) A key feature of these
bodies is that the {\it only} way to tile the plane with them is with
a tiling whose symmetry group does not have a fundamental domain of
finite volume; this situation is the defining characteristic of
``aperiodicity'', and has led to renewed study of the symmetry of
tilings (and thus packings); see [Radi]. There are other significant
symmetry features of this example which will be discussed below, where
they will help to develop an appropriate notion of equivalence among
optimally dense packings.

In practice it is almost impossible to actually determine optimally
dense packings -- for instance, there is nothing yet proven qualitatively
of the densest packings in the Euclidean plane by regular pentagons of
fixed size (see however [KuKu]) -- and we follow the lead of Hilbert and others in
concentrating on general features of optima, such as their geometric
symmetries. (It is because of their qualitative symmetry features that
we attributed high significance to the aperiodic tilings.) We will
discuss below some optimization results for packings in hyperbolic
spaces, by balls in $\H^n$, and by a certain polygon in $\H^2$.

We note that, for optimization
in a Euclidean space, there is no difficulty of {\it existence} of an
optimum {\it density}, even though, since the limit of \eqdens\ does
not exist for some packings, we are not able to make a comparison
among {\it all} packings. One way to understand this is to use the
fact that the relative densities $D_{B_r(p)}(P)$ of {\it any} packing
$P$ can be well approximated by packings $P'$ that have compact
fundamental domains; intuitively, the supremum of the densities of
such symmetric packings $P'$ is the desired optimum density, and can
be shown to be achieved, {\it in the sense of} \eqdens\
\hs-.05 , by some packing which is a limit of such
symmetric packings. (See [BoR2] for a complete argument.)

The situation for hyperbolic space packings is much more complicated,
due to the fact that the volume of a ball of radius $R$ grows
exponentially in $R$ ([Bear], [Kato]). This has the consequence that a significant
fraction of the volume is near the surface of the ball, so it is by no
means clear that one could make a useful approximation using a packing
with cofinite symmetry group. Once one is prevented from reducing
the problem to such symmetric packings, one is confronted with the
difficulty of showing the {\it existence} of a limit such as \eqdens\
\hs-.02 for a purported optimal packing; see [BoR2] for a history of
this difficulty.  In summary, lack of proof of the existence of
appropriate limiting densities was an impediment to progress in the
study of optimal density (and, {\it a fortiori}, optimally dense
packings) in hyperbolic spaces for many years [BoR2]. This led to a
search for alternatives to the notion of optimal density. The best
known of these are that of ``solid'' packings, and ``completely
saturated'' packings [FeKu]. Both notions are defined through local
properties of packings. The quest for a local approach/alternative to
optimal density is perhaps reasonable given that the only practical
way yet devised to prove that a packing is optimal is to prove that it
is locally optimal in all local regions (as noted above for disk
packing in $\E^2$, and as the basis for the methods for ball packings
in $\E^3$ [Laga]). However, no local alternative has proven
satisfactory [FeKu], [Bowe].

So far we have concentrated on the question of {\it existence} of
solutions (i.e., optimally dense packings) for our general optimization
problem, especially the difficulties for packings in a hyperbolic
space. This existence problem was solved recently ([BoR1], [BoR2]),
where the main obstacle, the existence of limiting densities, was
obtained through an ergodic theory formalism outlined in the next
section.  One goal of this paper is analysis of the {\it uniqueness}
of solutions for our optimization problems. Currently, there are
difficulties even for Euclidean problems. For instance, even though it
would be intuitively satisfying to declare that the problem of
optimally dense packings of $\E^2$ by disks of fixed radius has the
``unique'' solution discussed above (Figure 1), there has been no
satisfactory way to exclude some other packings of the same density,
for instance those obtained by deleting a finite number of disks from
this packing [CoGS]. This has been a serious obstacle to treatment of
optimal density as one treats other optimization problems [Kupe], and
will be a useful guide for our approach to a general
understanding of the qualitative features of optimally dense packings.

The situation is significantly more complicated, and interesting, for
the optimal packings (tilings) by kites \& darts in $\E^2$ (Figure 2)
than it is for equal disks. It can be proven that there are
uncountably many pairwise noncongruent such tilings, but that every
finite region in any one tiling appears in every other such tiling, so
they are, in some sense, ``locally indistinguishable'' [Gard]. It is
natural to want to declare that this optimization problem also has a
``unique'' solution, and this has, in effect, been the practice of
those studying aperiodicity, a practice we will follow.

The notion of uniqueness must differentiate between the situations for
densest packings of equal balls in $\E^3$ and that of $\E^2$; for the
former the expected solution is intuitively far from unique,
containing for instance the face centered cubic and also the hexagonal
close packed structures, which must be considered different if the
notion is used at all. However at a deeper level it must also give a
useful criterion for when two such optimal structures are ``the
same'', that is, it must give some useful notion of the geometric
symmetry of optimally dense packings, for instance of the kite \& dart
tilings. In a later section we will discuss the connection between our
approach and the use, by Connes and others, of noncommutative topology
to understand the symmetry of structures such as the kite \& dart
tilings; see the expository works [Conn], [KePu] and references they 
contain.

As we hope to demonstrate, study of the uniqueness problem for
packings of hyperbolic space will be useful even for understanding
packings in Euclidean space. Much of this paper consists of the
analysis of a specific family of examples of optimal density in a
hyperbolic space, the {first} aperiodic examples in hyperbolic space
for which explicit
optimally dense packings have been determined. (See [BoR2] for
an analysis of other examples of tilings in $\H^2$, such as those in
[MaMo], [Moze] and [Good]). These examples exhibit
features not seen in Euclidean examples, and will help us draw some
conclusions about the general features of optimum density
problems.
\vs.2

\nd {\bf II. The Ergodic Theory Formalism for Optimal Packing}
\vs.1
One concern of this paper is with the uniqueness of solutions to optimal
density packing problems. While proving existence of solutions for
problems set in Euclidean spaces did not require giving the problems a
formal structure, the question of existence of solutions for problems
set in hyperbolic space definitely did require introducing a formal
structure, and we will see that this same structure is useful for
handling questions of uniqueness, even in Euclidean space. We follow
[BoR1], [BoR2] in introducing an ergodic theory structure into our
optimization problems, in order to control the existence of limits
such as \eqdens\ \hs-.05 .

Using the notation of section I, consider the space $\P_\B$ of all
possible packings of $\S^n$ by bodies from $\B$, and put a metric
on $\P_\B$ such that convergence of a sequence of packings
corresponds to uniform convergence on compact subsets of $\S^n$. Such
a metric makes $\P_\B$ compact, and makes continuous the natural
action on $\P_\B$ of the (connected) group $\G^n$ of rigid motions of
$\S^n$ [RaWo].

We next consider Borel probability measures on $\P_\B$ which are
invariant under $\G^n$. (To construct such a measure consider any
packing $P$ for which the symmetry group has fundamental domain of
finite volume, and identify the orbit $O(P)$ of $P$ under $\G^n$ as
the quotient of $\G^n$ by that symmetry group. One can then project
Haar measure from $\G^n$ to an invariant probability measure on $O(P)$
and then extend it to all of $P_\B$ so that the complement of $O(P)$
has measure zero.)

We define the ``density of the (invariant) measure $\mu$''
on $\P_\B$, $D(\mu)$, by $D(\mu)\equiv\mu(A)$, where $A$ is the following set of
packings:
$$A\equiv \{P\in \P_\B\,|\, \hbox{the origin } \O \hbox{ of } \S^n
\hbox{ is in a body in }P\}. \cheqno\eqA $$

\nd (It is easy to see from the invariance of $\mu$ that $\mu(A)$ 
is independent of the choice of origin.) We may now introduce the
notion of optimal density.
\vs.1 \nd
{\bf Definition 1}. A probability measure $\bar \mu$ on the space
$\P_\B$ of packings, ergodic under rigid motions, is ``optimally
dense'' if $D(\bar \mu)=\sup_{\mu} D(\mu)=\sup_{\mu}\mu(A)$; the
number $\sup_{\mu}\mu(A)$ is the ``optimal density'' for packing
bodies from $\B$.
\vs.1 \nd
(An invariant measure $\mu$ is ``ergodic'' if it cannot be expressed
as an average: \hbox{$\mu=a_1\mu_1 +a_2\mu_2 $,} with $a_1,\, a_2>0$
and $\mu_1 \ne \mu_2$ invariant.) It is not hard to show
[BoR1] the existence of such optimal measures for any given $\B$.

Finally, the terminology is justified as follows. First we rewrite the
right hand side of \eqdensrho\ as:

$$\lim_{r\to \infty}{1\over \nu[\G^n(r,p)]} \int_{\G^n(r,p)} \chix_A[g(P)]\, d\nu(g), 
\cheqno\eqquot $$

\nd where $\chix_A$ is the indicator function for $A$, $\nu$ is
Haar measure on $\G^n$ and \hfill
$$\G^n(r,p) = \{g\in \G^n\,| \, d_{\S^n}[g(p),p]<r \}, \cheqno\eqGd $$

\nd where $d_{\S^n}$ is the distance function on $\S^n$. It follows from
G.D. Birkhoff's pointwise ergodic theorem [Walt] that for $\S^n=\E^n$ and any ergodic
$\mu$ there is a set of $P$'s, of full $\mu$-measure and invariant
under $\G^n$, for which the limit in \eqquot\ exists. This has been
extended to $\S^n=\H^n$ by Nevo et al.: [Nevo], [NeSt] (with the
invariance of the set of $P$'s proven in [BoR2]).  We may conclude
then that ``most'' of the packings in the support of a fixed ergodic
measure have the same well defined density in the sense of \eqdens\
\hs-.05 ; so as one varies the
measure one sees $\P_\B$ decomposed into packings of various
densities, with those of optimal density being the ones in which we
are interested.  Formally we define optimally dense packings (slightly
more stringently than in [BoR2]) as follows.
\vs.1 \nd
{\bf Definition 2}. A packing $P$ is ``optimally dense'' if it is 
``generic'' for some optimally dense $\mu$, that is, it is in 
the support of $\mu$ and:
$$ \int_{\P_\B} f(Q)\,d\mu(Q)= \lim_{r\to\infty} {1\over
\nu[\G^n(r,p)]} \int_{\G^n(r,p)} f[g(P)]\, d\nu(g), \cheqno\eqoptpack$$

\nd for every $p \in \S^n$ and every continuous $f$ on $\P_\B$. The set
of all optimally dense packings for bodies in $\B$ will be 
denoted $\P^o_{\B}$.
\vs.1 
We note that the set of packings generic for the invariant measure $\mu$ is
of full measure with respect to $\mu$; this follows from the ergodic theorem
and the fact that the space of continuous functions on $\P_\B$ is separable
in the uniform norm.
\vs.1 \nd
{\bf Lemma 1}. If $P$ is generic for the ergodic measure $\mu$ (which
is not necessarily optimally dense) then:
$$ \eqalign{\lim_{r \to\infty} {1\over
\nu[\G^n(r,p)]} \int_{\G^n(r,p)} \chix_A[g(P)]\, d\nu(g)&= 
\int_{\P_\B} \chix_A(Q)\,d\mu(Q)\cr
&= D(\mu)\cr}
\cheqno\eqgeneric$$

\nd for every $p\in \S^n$.
\vs.1 \nd
Proof. Let 
$$A'\equiv \{P \in \P_\B\,|\,  \O 
\hbox{ is in the interior of a body in }P\}. \cheqno\eqap $$

\nd Define continuous $f_k$ on $\P_\B$ by:

$$f_k(P)=\cases{1, \ \ \hbox{on } A \cr
0, \ \ \hbox{on } \{P \in A^c\,|\,d_{\S^n}(\O,\partial P)\ge {1\over k}\}\cr
1-kc,\ \ \hbox{on } \{P\in A^c\,|\,d_{\S^n}(\O,\partial P)=c < {1\over k}\},\cr} 
\cheqno\eqab $$
\nd where $\partial P$ denotes the union of the boundaries of the bodies in $P$.
Note that the $f_k$ decrease pointwise to $\chix_A$. Similarly define 
continuous $g_k$ on $\P_\B$ by:

$$g_k(P)=\cases{0, \ \ \hbox{on } {A'}^c \cr
1, \ \ \hbox{on } \{P\in A'\,|\,d_{\S^n}(\O,\partial P)\ge {1\over k}\}\cr
kc,\ \ \hbox{on } \{P\in A'\,|\,d_{\S^n}(\O,\partial P)=c < {1\over k}\}.\cr} 
 \cheqno\eqac $$
\nd
Note that the $g_k$ increase pointwise to $\chix_{A'}$. Now,
given $\epsilon >0$, choose $K > 0$ such that 

$$0<\int f_K \,d\mu - \int \chix_A \,d\mu < \epsilon/2 \ \ \hbox{  and  } \ \ 
0<\int \chix_{A'} \,d\mu - \int g_K \,d\mu < \epsilon/2. \cheqno\eqad $$

\nd Define, for $R>0$, the measure $\nu(R,P)$ on $\P_\B$ by
$$\int f \,d\nu(R,P) \equiv {1\over \nu[\G^n(R,\O)]} \int_{\G^n(R,\O)} f[g(P)]\ d\nu(g), 
\cheqno\eqae $$
\nd for continuous $f$, where as before $\nu$ is Haar measure on $\G^n$. Then choose 
$\tilde R>0$ such that 
$$\Big|\int f_K \,d\nu(R,P) - \int f_K \,d\mu\Big|< \epsilon/2 \ \ \hbox{  and  } \ \ 
\Big|\int g_K \,d\nu(R,P) - \int g_K \,d\mu\Big|< \epsilon/2 
\cheqno\eqaf $$

\nd for all $ R > \tilde R  $. We then have:
$$\eqalign{\int \chix_{A'} \,d\mu - \epsilon < \int g_K \,d\nu(R,P) <
\int &\chix_{A'} \,d\nu(R,P) 
\le \int \chix_A \,d\nu(R,P) \cr
< \int f_K \,&d\nu(R,P) < \int \chix_A \,d\mu + \epsilon.\cr} \cheqno\eqag $$

\nd However, from the ergodic theorem $\int \chix_A-\chix_{A'} \,d\mu =
\int \chix_{A/A'}\,d\mu=0$, 
so $\int \chix_{A'} \,d\nu(R,P)$ and $\int \chix_A \,d\nu(R,P)$ both
converge to $\int \chix_A \,d\mu=D(\mu)$ as $R\to \infty$.\qed

\vs.1 
Next we note a useful tool for computing optimal densities.
For those $P\in \P_\B$ such that the point $p\in S^n$ is contained in
the interior of a Voronoi cell (which cell we denote by $V_p(P)$), we
define $F_p(P)$ to be the relative volume of $V_p(P)$ occupied
by the bodies of $P$. (We note that $F_p$ is defined
$\mu$-almost everywhere for any invariant $\mu$.)
\vs.1 \nd
{\bf Definition 3}. For invariant measures $\mu$ we define 
the ``average Voronoi density
for $\mu$'', $D_V(\mu)$, as $\int_{\P_\B} F_p(P)\,d\mu(P)$.
(Note that $D_V(\mu)$ does not depend on $p$ because of the invariance of $\mu$.)
\vs.1 
The notion of average Voronoi density is useful, as it has been shown
[BoR1, BoR2] that, for any invariant measure $\mu$, the average Voronoi
density $D_V(\mu)$ equals the average density $D(\mu)$.
\vs.1 
Summarizing the above, we have sketched a formalism through which one
proves {\it existence} of solutions to the general problem of densest
packings of $\S^n$ by congruent copies of bodies from $\B$. Our 
next goal is to consider uniqueness, but we will first need to discuss
symmetry further (section III), and a new family of examples (section IV).
\vs.4 \nd
{\bf III. Symmetry in the Problem of Optimally Dense Packings}
\vs.1
We will use the following common terms: a packing is called
``periodic'' (resp. ``nonperiodic'') if it has (resp. does not have) a
symmetry group with fundamental domain of finite volume, and we say an
optimal packing problem is ``aperiodic'' if {\it all} its optimally
dense packings are nonperiodic.

Although it will be more difficult to introduce a general criterion for
uniqueness (which we will attempt in section V), the above
formalism can be used to solve the old problem of making sense of
uniqueness for the problem for disks of fixed radius in $\E^2$. (One
consequence of this difficulty was the search for replacements of
the notion of density, such as ``solidity'' and ``complete
saturation''; see [FeKu], [Feje], [CoGS] and [Bowe].)
\vs.1 \nd
{\bf Theorem 1}. There is only one optimally dense packing in $\E^2$
for disks of fixed radius, up to rigid motion.
\vs.1 \nd
Proof. Assume $\mu$ is any optimally dense measure for this
problem. Using the fact that the hexagonal Voronoi cells of $ P^o$
(the hexagonal packing of Figure 1) are, up to rigid motion, the unique cells of
optimal density (see [Feje]), and using the basic result on average Voronoi cells
[BoR1], we see that for $\mu$-almost every packing $P$, the cell
containing a particular point $p$ must be this regular hexagon. Repeating this
argument for a countable dense set of $p$'s, we see that for
$\mu$-almost every packing $P$ {\it every} Voronoi cell is, up to
rigid motion, this regular hexagon, i.e., $\mu$-almost every packing
$P$ is $P^o$, up to rigid motion. It is shown in [BoR1] that
there is a unique invariant measure with support in the (closed) orbit
of any periodic packing, and therefore the orbit of $P^o$
consists of all the optimal packings.\qed
\vs.1
The ergodic theory formalism automatically makes sense of the
uniqueness of this optimization problem, by its subjugation of
packings to invariant measures on packings. Aperiodic problems, such
as the kites \& darts, are more subtle. Our next result is a
connection between aperiodicity and the uniqueness of packing
problems.
\vs.1 \nd
{\bf Theorem 2}. If there is only one optimally dense packing of
$\S^n$, up to congruence, by congruent copies of bodies from some
fixed, finite collection $\B$, then that packing must have a symmetry
group with compact fundamental domain.
\vs.1 \nd 
Proof. By assumption there exists a probability measure $\mu$,
invariant under $\G^n$, for which the orbit $O(P)$ of some packing $P$
has measure one. Since $O(P)$ can be identified with the quotient of
$\G^n$ by the symmetry group $\Gamma_P$ of $P$, it follows from
the uniqueness of Haar measure on $\G^n$ that
$\Gamma_P$ is cofinite.
 
Since the volume of $\S^n/\Gamma_P$ is finite, only a finite
number of bodies can appear in any fundamental domain, and in particular
the bodies lie in a compact region of $\S^n/\Gamma_P$.  
If $\Gamma_P$ were cofinite but not cocompact (something possible
only if $\S^n$ is hyperbolic), then the preimage in $\S^n$ of the ends of
$\S^n/\Gamma_P$ would not contain any bodies. 
However, the preimage of a hyperbolic end contains arbitrarily large
balls, so there is room in our packing $P$ to add additional bodies.

This contradicts the fact that optimally dense packings must be 
``saturated'', meaning that 
one cannot add another body to the uncovered regions.
In fact it was proven in [Bowe] that the set of
{\it completely} saturated packings have full measure with respect to
any optimally dense measure. (A ``completely saturated packing'' is one in
which it is impossible to remove any finite number of bodies 
and replace them with bodies of larger total
volume.) \qed
\vs.1 
One thing we can conclude from these two theorems is that the problems of
optimal density decompose naturally into two classes: those
allowing periodic optima (such as with balls of fixed size in
$\E^2$ or $\E^3$), and the class of aperiodic problems. The former no
longer pose any difficulty as to classifying their uniqueness, leaving
us now to understand the more interesting class of aperiodic problems.

Aperiodicity is not an unnatural circumstance -- it may even be
generic in some sense; see [MiRa] for a related problem. In Euclidean
space aperiodicity has only been discovered so far in packings of
complicated polyhedra, whereas in hyperbolic space it already appears
in packings of balls.
\vs.1 \nd
{\bf Theorem 3 [BoR1]}. For all but countably many fixed radii $R$ the
ball packing problem in $\H^n$ has only aperiodic solutions.
\vs.1 
For packings/tilings in Euclidean space there has not yet been a
significant attempt to understand the uniqueness problem. We will
postpone our attempt at a formal definition until after considering
the following new examples, the first aperiodic problem in hyperbolic
space for which we can determine explicit solutions. 

\vs.2 \nd
{\bf IV. A New Example: the Modified Binary in $\H^2$}
\vs.1
Using the upper half plane model of $\H^2$, consider the ``binary
tile'' $\tau$ of Figure 4 (introduced by Roger Penrose in 1978
[Penr]).  We define as the ``core'' $\gamma$ of $\tau$ that shape with
four edges -- two segments of geodesics, and two segments of
horocycles (one twice the length of the other) -- obtained by omitting
the bumps and dents of $\tau$. More specifically we could take the
coordinates of the vertices of $\gamma$ to be $i,\ i+2, 2i$, and
$2i+2$. Congruent copies of $\tau$ can only tile $\H^2$ as in Figure
5. In fact it is useful to classify these tilings as follows. Once the
location of one tile is known, the bumps on the geodesic edges of (the
core of) $\tau$ force the positions of tiles filling out the region
between two concentric horocycles (the ones containing the horocyclic
edges of the (core of the) tile). Consider now the possible tiles
abutting the ones in this ``horocyclic strip''. There is only one way
to fill an abutting strip which is ``further'' from the common point
at infinity of the horocycles, and two ways to fill the strip which is
``closer''. This fully classifies the possible tilings of $\H^2$ by
$\tau$. 

We now construct a new tile $\bar \tau$ based on the same core $\gamma$, which
will permit some new tilings. On the geodesic edges of $\gamma$ we add the
same bumps and dents as before, but we enlarge each of the other
original bumps and dents as follows. Consider a (densest [Boro])
packing of the hyperbolic plane by horoballs, as illustrated in Figure
6, and consider three abutting horoballs. Divide the (white) region
between the horoballs into 3 congruent regions by means of three
geodesics, each drawn from the center of three-fold symmetry of the region to the
point where a pair of horoballs touch. The new tile $\bar \tau$ is
illustrated in Figure 7.
 From the construction we see that, as with $\tau$, in any tiling of the
plane by $\bar \tau$, once we know the location of a specific tile we can
uniquely fill in a horocyclic strip, and then have two choices for
filling, consecutively, each of those strips which are closer to their
point at infinity, thus filling in a horoball. However,
there is now a second way to fill the abutting strip which is further
from the point at infinity, in 
which the bumps from three tiles abut to fill in a white region in Figure 6. If
this latter method is used, the only way to complete a tiling of the
plane is to fill in each of the horoballs defined by the two new tiles, 
then use the same method to extend beyond these horoballs
to more horoballs, etc. Intuitively, these new tilings are obtained from a densest horoball 
packing by tiling each horoball with copies of the tile ${\bar \tau}$.

Let $\T({\bar \tau})$ be the set of all tilings of the plane by $\bar
\tau$.  We call such a tiling ``degenerate'' if the cores of the tiles
themselves tile the hyperbolic plane, and ``non-degenerate'' if the
union of the cores corresponds to a (densest) packing $\bar P$ of
$\H^2$ by horoballs.  (Such a horoball packing has symmetry group conjugate to
PSL$(2,\Z)$, 
which is cofinite but not cocompact.) 
We know from [BoR1] that the set $\T_{deg}$ of
degenerate tilings has measure 0 with respect to
any invariant probability measure on $\T({\bar \tau})$.  As a result, 
degenerate tilings do not qualify as optimally dense packings of $\bar \tau$
in the sense of Definition 2.  However, invariant measures on $\T({\bar \tau})$
do exist.

To construct such a measure we consider the internal structures
of the different horoballs. This structure is related to a
choice of the two lower bumps on the tile $\bar \tau$, henceforth
called ``prongs''.  For each horoball $H$ and each triangle that
touches the horoball, the internal structure of $H$ is associated to a
dyadic integer, that is a formal sum $\sum_{i=0}^\infty a_i 2^i$, 
with $a_i\in\{0,1\}$, where two dyadic integers are
considered close if their first $N$ terms agree, with $N$ large. 
The first digit tells whether the left or
right prong of a tile from $H$ sticks into the triangle, the next digit tells 
whether that tile
emerges from the left or right prong of its ``parent'', and so on.
(We will eventually use the algebraic structure of these quantities.) 
We let each digit be an independent random variable, with equal
probability of being 0 or 1. Let the internal structures of distinct
horoballs be independent, and be independent of the location of the 
horoballs, which is given by Lebesgue measure on $\H^2/$PSL$(2,\Z)$.

 From the existence of an invariant measure, it follows [BoR1] that 
there are tilings in $\T({\bar \tau})$ which are
optimally dense packings of $\bar \tau$ in the sense of Definition
2. It is not difficult to show from this that {\it all} optimally
dense packings of $\bar \tau$ are tilings [BoR1]. 
We now see that the optimal density
problem for the modified binary $\bar \tau$ is aperiodic: the presence
of the horoballs immediately implies that the symmetry group of an
optimal packing/tiling is at most cofinite, not cocompact, and, as
argued in the proof of Theorem 2, a tiling by compact bodies cannot
have a symmetry group which is cofinite but not cocompact.

The presence of the closed invariant set $\T_{deg}$ is a new feature
in optimal density problems. Since it is a subset of the orbit closure
of every tiling, it is in the support of every invariant measure on
$\T({\bar \tau})$, even the ergodic ones. However, there is no invariant
measure on $\T_{deg}$. This is not possible for problems set in a
Euclidean space, since if the Euclidean group acts on a compact metric
space an elementary fixed point argument ([Radi]) guarantees the
existence of an invariant probability measure on that set. We will
discuss this feature of $\T_{deg}$ further below, when we consider
various types of conjugacy for the dynamical systems in which we are
couching our optimization problems.

We now turn to the construction of uniquely ergodic invariant subsets of 
$\T({\bar \tau})$, and the measures they support. 
Let $w$ be a dyadic integer.
For each $w$, let $\T_w$ be the closure of the class of
tilings for which the sum of the three dyadics at each triangle
is $w$.

\vs.1 \nd
{\bf Theorem 4}. $\T_w$ is uniquely ergodic under the action
of $\G^2=\,$PSL$(2,\R)$.
\vs.1 \nd
Proof.  For $T \in \T_w - \T_{deg}$, there is naturally associated to
$T$ a horoball packing $h(T)$. For $N>0$ and $T \in \T_w - \T_{deg}$,
we let $\Theta_N(T)$ be the packing obtained from $T$ by removing all
but the $N$ 
horocyclic rows of tiles closest to the
boundary of any horoball in $h(T)$. If $T \in \T_{deg}$, we let
$\Theta_N(T)$ be the empty packing, $\emptyset$. Note that $\Theta_N$ defines a
continuous map from $\T_w$ onto a compact space $\P_N$ of packings and
that $\Theta_N$ commutes with the action of $\G^2$.  

 From Lemma 2 below it follows that $\P_N$ admits only two ergodic measures; 
one concentrated on $\emptyset$
and the other being derived from Haar measure on the space 
$\G^n/H$ where $H$ is the symmetry group of $\Theta_N(T)$ for any $T \in 
\T_w - \T_{deg}$ ([BoR2]). Since $\T_{deg}$ has $\mu$-measure zero with respect
to any invariant measure $\mu$ on $\T_w$ ([BoR2]), 
the empty packing is measure zero with respect to the pushforward of 
$\mu$.
Therefore there is only one possibility for the 
pushforward of $\mu$, or equivalently, the space $\T_w$ is
uniquely ergodic with respect to the $\sigma$-algebra
$\Theta_N^{-1}(\Sigma_N)$ (where $\Sigma_N$ is the Borel
$\sigma$-algebra of $\P_N$). Since this is true for all $N$ and the
$\sigma$-algebras $\Theta_N^{-1}(\Sigma_N)$ are increasing, this
implies that $\T_w$ is uniquely ergodic with respect to the
$\sigma$-algebra $\bigcup_N \Theta_N^{-1}(\Sigma_N)$. But $\T_{deg}$
has measure zero with respect to any invariant measure $\mu$ on $\T_w$
([BoR2]) and the topology of $\T_w - \T_{deg}$ is generated by $\cup_N
\, \Theta_N^{-1}(\t_N)$ (where $\t_N$ is the topology on $\P_N$). So
the $\mu$-closure of $\bigcup_N \, \Theta_N^{-1}(\Sigma_N)$ contains
the Borel $\sigma$-algebra of $\T_w$. Thus $\T_w$ is uniquely
ergodic. It remains only to prove the following lemma.
\vs.1 \nd
{\bf Lemma 2}. For any $T,\, T' \in \T_w -\T_{deg}$, $\Theta_N(T)$ has
cofinite symmetry group, $\Theta_N(T)$ is in the orbit of 
$\Theta_N(T')$ and $\P_N$ is equal to this orbit union
the empty packing.
\vs.1 \nd
Proof. The symmetry group of the horoball packing $h(T)$
for $T \in \T_w-\T_{deg}$ is conjugate to PSL$(2,\Z)$.  The packing
$h(T)$ naturally corresponds to an infinite trivalent tree with the
triangles of the packing corresponding to the vertices of the tree. We
can navigate around the tree with two fundamental operations.

 If a ``state'' is a vertex together with a
choice of one of the three edges leading out from that vertex, then
the two operations on states are

$$ C = \hbox{Rotate counterclockwise by 120 degrees} \cheqno\eqprfa $$
$$ L = \hbox{Go forwards to the next vertex and bear left}. \cheqno\eqprfb $$

\nd $C$ and $L$ obviously generate the entire symmetry group of the tree.
In terms of PSL$(2,\Z)$, $C$ is the elliptic element 
$$C=\pmatrix{0 & 1\cr -1 & 1} \cheqno\eqprfc$$ 

\nd or $z \to 1/(1-z)$, and $L$ is the parabolic element
$$L=\pmatrix{1 & 1 \cr 0 & 1} \cheqno\eqprfd $$ 
\nd 
or $z \to z+1$, and together they generate
all of PSL$(2,\Z)$.  

Now we consider filling the horoballs with tiles $\bar \tau$, so
that to each triangle we can associate three dyadic integers, one for
each of the horoballs that meet at the triangle.  Different triangles
that touch the same horoball will not have the same dyadic integer;
rather, moving along the edge of the horoball counterclockwise (as
seen from inside the horoball) will increase the dyadic number by one
each step.

If we apply the condition that the three numbers at each triangle must
add up to $w$ (a fixed dyadic integer), then two dyadic numbers at a
triangle determine all the rest.  If you know that a given vertex has
numbers $a$ and $b$ then the third number must be
$c=w -a-b$. One neighboring vertex has numbers $b-1$ and $a+1$, so its
third vertex must be $c$.  Each vertex determines its neighbors,
and so determines the entire tree.  Thus a tiling in standard position
(that is, with a choice of preferred vertex and preferred edge
directed out from the preferred vertex) can be associated to a pair
$(a,b)$ of dyadic numbers, where $a$ is the index of the horoball to
the left of the preferred outgoing edge, and $b$ is the index of the
horoball to the right.

The action of $L$ and $C$ is easy to compute, namely:
$$L : (a,b) \mapsto (a+1, c) \cheqno\eqprfe$$
$$C : (a,b) \mapsto (c, a), \cheqno\eqprff$$
\nd where $c = w -a-b$.   It is not hard to check the following
elements:
$$ L^2 : (a,b) \mapsto (a+2, b-1) \cheqno\eqprfg$$
$$ C L^2 C^2 : (a,b) \mapsto (a-1, b+2) \cheqno\eqprfh$$
$$ L^4 C L^2 C^2: (a,b) \mapsto (a+3, b) \cheqno\eqprfi$$
$$ L^2 C L^4 C^2: (a,b) \mapsto (a, b+3). \cheqno\eqprfj$$

Now consider the effect of PSL$(2, \Z)$ on pairs $(a,b)$ as above but
taken modulo $2^N$. Since 3 and $2^N$ are relatively prime, some power
of $L^4 C L^2 C^2$ sends $(a,b)$ to $(a+1,b)$\hs-.05 $\pmod{2^N}$, and
some power of $L^2 C L^4 C^2$ sends $(a,b)$ to $(a,b+1)$\hs-.05
$\pmod{2^N}$.  Thus PSL$(2,\Z)$ acts transitively on the space of
pairs $(a,b)$\hs-.05 $\pmod{2^N}$.

Therefore, for all $T,T' \in \T_w -\T_{deg}$, $\Theta_N(T)$ is
congruent to $\Theta_N(T')$. Also, the subgroup that preserves the
pair $(a,b)$\hs-.05 $\pmod{2^N}$ (and hence the first $N$ rows of each
horoball in the tiling) is an index $2^{2N}$ subgroup of PSL$(2,\Z)$,
and hence is a cofinite subgroup of PSL$(2,\R)$. Finally, $P_N$
is the union of this orbit and the image of $\T_{deg}$, which is
the empty packing.\qed

\vs.1 

We next give another property of these ``fixed-sum'' classes of
tilings, in terms of dynamical conjugacy. For convenience we recall
some common terms. A topological group $G$ acts continuously on a
compact metric space $X$ if there is a map $\phi: (g,x)\in G\times
X\to g(x)\in X$ which is continuous and satisfies $h[g(x)]=[hg](x)$
for all $g,h\in G$ and $x\in X$. Assuming $G$ acts continuously on $X$
and $Y$, the actions are called ``topologically conjugate'' if there
is a homeomorphism $\alpha: x\in X\to \alpha[x]\in Y$ such that
$\alpha [g(x)]=g(\alpha[x])$. Assume further the existence on $X$ and
$Y$ of Borel probability measures $\mu_X$ and $\mu_Y$ which are
invariant under the corresponding actions of $G$. These two actions
are said ``measurably (or metrically) conjugate'' if there are
invariant subsets $X_0\subset X$ and $Y_0\subset Y$, each of measure
zero, and an invertible map $\alpha': x\in X/X_0\to \alpha'[x]\in
Y/Y_0$ such that $\alpha'[g(x)]=g(\alpha'[x])$ which, together with
its inverse, is measure preserving. Finally we introduce an intermediate
form of conjugacy (related to ``almost topological conjugacy'' [AdMa])
as follows. The actions of $G$ on $X$ and $Y$ will be called ``almost
conjugate'' if there are
invariant subsets $X_0\subset X$ and $Y_0\subset Y$, each of measure
zero with respect to all invariant measures, 
and a homeomorphism $\alpha': x\in X/X_0\to \alpha'[x]\in
Y/Y_0$ such that $\alpha'[g(x)]=g(\alpha'[x])$.
\vs.1 \noindent 
{\bf Theorem 5}. $\T_w$ and $\T_{w'}$ are topologically conjugate
if and only if $w-w' \in 3\Z$. $\T_w$ and $\T_{w'}$ are almost
conjugate for any $w,w'$.
\vs.1 \nd
Proof. If $w'-w =3e$, where $e \in
\Z$, then we construct a conjugacy by leaving the location of all 
the horoballs fixed, and simply adding $e$ to the dyadic index of 
each horoball.
In the $N$-th layer of a horoball the conjugacy essentially acts by
translation by $e/2^N$, so points deep within a horoball are moved
only slightly.  In the case of degenerate tilings, the conjugacy
leaves the entire tiling fixed.

If $w-w'$ is not 3 times a rational integer, then $(w-w')/3$ is still a dyadic 
integer, since 3 is a unit in the ring of dyadic integers. Adding $(w-w')/3$
to the index of each horoball is a continuous map on the complement
of $\T_{deg}$, but is not uniformly continuous and does not extend to
all of $\T_w$. This shows that $\T_w$ and $\T_{w'}$ are almost conjugate
for any $w,w'$. 

The proof that $w-w' \in 3\Z$ is necessary for topological conjugacy
is harder, and will consist of four lemmas.
\vs.1 \nd
{\bf Lemma 3}. Any topological conjugacy $\phi$ between $\T_w$ and
$\T_{w'}$ must preserve the points on the sphere at infinity that
are tangent to horoballs.  Furthermore, the ``radii'' of the horoballs
can only change by a finite amount.  That is, there exists a constant
$R$ (depending only on $\phi$) such that, if $T\in \T_w$ is a tiling and 
$H$ is a horoball in $h(T)$,  and $H'$ is the corresponding horoball in
$h[\phi(T)]$ (that is, with the same tangent point), 
then $H$ is contained in an $R$-neighborhood of $H'$ and
vice-versa.
\vs.1 \nd
{\bf Lemma 4}. The $R$ of the previous lemma is actually zero;
topological conjugacies preserve the locations of horoballs exactly.
\vs.1 \nd
{\bf Lemma 5}. Let $H$ be any horoball in $h(T)$ for any $T \in \T_w$, and
let $a$ be its dyadic index (measured from a particular triangle). Let
$a'$ be the dyadic index of $H' \subset h[\phi(T)]$ measured from the
same triangle. Then the set of differences $a'-a$ (for all such
horoballs $H$ in all such tilings $T$) is a bounded subset of $\Z$.
\vs.1 \nd
{\bf Lemma 6}. There is a triangle in $T$, with indices $a$, $b$ and $c$, such
that $a'-a = b'-b = c'-c = (w'-w)/3$.
\vs.1 \nd
Proof of Lemma 3. A topological conjugacy is uniformly continuous, so
if $\phi: \T_w \to \T_{w'}$ is a topological
conjugacy, then for every $r'>0$ and $\epsilon>0$ there is a radius
$r$ such that, if the neighborhoods of two points agree to radius $r$,
then from uniform continuity and
conjugacy their images agree out to radius $r'$, up to an ``$\epsilon$
wiggle''. More precisely, for any two tilings $T_1, T_2 \in
\T_w$ and points $p_1$, $p_2\in \H^2$, if there is an isometry of $\H^2$
that sends a ball of radius $r$ of $p_1$ in $T_1$ exactly onto a ball
of radius $r$ of $p_2$ in $T_2$, then the same isometry sends a ball
of radius $r'$ of $p_1$ in $\phi(T_1)$ to an $\epsilon$-small
distortion of a ball of radius $r$ of $p_2$ in $\phi(T_2)$, where an
``$\epsilon$-small distortion'' means an isometry that moves each
point in the neighborhood a distance $\epsilon$ or less.

Now take $r'$ to be greater than the diameter of a triangle and
$\epsilon$ to be much less than the diameter of a triangle. Take any
tiling $T$ for which $\phi(T)$ has a triangle centered at $p_2$. We
claim that $T$ contains a triangle centered at a point $p_1$ at
distance at most $r+1$ from $p_2$.  For if not, then the $r$-neighborhood of 
$p_2$ (call it $U$) lies 
completely within a horoball of $h(T)$. But then there is a constant, 
${\tilde r}$, say such that every ball of size ${\tilde r}$ contains an 
$r$-ball such that $T$ restricted to that $r$-ball is isometrically 
conjugate to $T$ restricted to $U$. (Here $T$ is thought of as the 
function from the plane to the tile ${\bar \tau}$ that is induced by the 
tiling $T$). This implies that every ball of size ${\tilde r} + 
\epsilon$  contains a triangle of $\phi(T)$. But this contradicts the 
fact that there are points, deep within a horoball of $h(\phi(T))$ that 
are at least a distance ${\tilde r} + \epsilon$ away from any triangle 
in the complement of $h(\phi(T))$. But this contradicts the fact that 
$\phi(T)$ is made up
of horoballs, some points of which are arbitrarily far from triangles.

Now let $H$ be any horoball in $h(T)$.  Since all points in $H$ that are
farther than $r+1$ from the boundary of $H$ are mapped into a horoball 
in $h[\phi(T)]$ (i.e., not into a triangle), and since this set of points 
is connected, $H$ lies within
an $r$-neighborhood of a specific horoball $H' \subset h[\phi(T)]$.  
This implies that $H$ and $H'$ have the same tangent point on the 
sphere at infinity.

To obtain the fixed bound $R$, just repeat the argument for $\phi^{-1}$ and
take $R$ to be the larger of the two constants $r+1$.\qed
\vs.1 \nd
Proof of Lemma 4.  Let $R$ be as before. We know that $\phi$
preserves the location of each horoball, and changes its radius by at
most $R$.  The question is which horoballs grow and which shrink, and
by how much.  Consider a triangle in $T$, with center point $p$, where
horoballs $H_1$, $H_2$ and $H_3$ meet. It is impossible for two of
these horoballs to grow, or one to grow while a second does not change,
lest they overlap.  If two stay fixed, then the entire pattern is
fixed.  Thus, if there are {\it any} changes {\it anywhere}, then at
each triangle either one horoball grows (or stays fixed) 
and the other two shrink, or
all three shrink and one or more other horoballs $H_4, H_5, \ldots$ grow 
to fill up the
space.  There are only a finite number of horoballs within distance
$R$ of the triangle, so only a finite number of possible directions
where the tile containing $p$ in $\phi(T)$ can be pointing. By 
continuity
there is a number $N$ such that knowing all tiles within the $N$th
collar of the triangle determines which of $H_1$, $H_2$ and $H_3$ grow
and which shrink.  In particular, knowing the first $N$ digits of the
dyadic labels for $H_{1,2,3}$ determines which grow and shrink.  In
essence, all our labels should be counted mod $2^N$.

There are $2^{2N}$ possible triples $(a,b,c)$ of numbers \hs-.15 $\pmod{2^N}$
that add up to $w$ \hs-.15 $\pmod{2^N}$. These correspond to labels for
horoballs that meet at a triangle, counting clockwise.  For each one,
either all triangles with this label have the ``$a$'' tile shrink, or
none do.  Let $S_a$ be the set of labels for which the ``$a$'' tile
shrinks, let $S_b$ be the set for which the ``$b$'' tile shrinks, and
let $S_c$ be the set for which the ``$c$'' tile shrinks.  We will show
all three sets are the whole set of triples, so all horoballs shrink,
which is a contradiction.

Each horoball actually meets an infinite number of triangles.  By
comparing adjacent triangles with the same horoball we see that
$(a,b,c)\in S_a$ if and only if $(a+1, c-1, b) \in S_a$. Continuing
this process, we get that either the entire orbit $\{ (a+2n, b-n,
c-n)\} \cup\{ (a+2m+1, c-m-1, b-m)\}$, is in $S_a$ or the entire orbit
is out.  Note that $3$ is a unit in $\Z_{2^N}$, so we can take
$n=(b+1-a)/3$ and $m=(b-a-2)/3$.  This means that both
$$\left ( {a+2b+2 \over 3}, {a + 2b -1 \over 3}, {c+a-b-1 \over
3}\right) \cheqno\prfba $$
\nd and
$$\left ( {a+2b-1 \over 3}, {c+a-b-1 \over 3}, {a +
2b +2 \over 3}\right) \cheqno\prfbb $$

\nd are in the same orbit as $(a,b,c)$.  If $(a,b,c) \not \in S_a$, then
in any triangle with indices ${(a+2b+2)/3}$, ${(a + 2b -1)/3}$, 
and ${(c+a-b-1)/3}$, in clockwise cyclic order, the horoballs
with indices ${(a+2b+2)/ 3}$ and ${(a + 2b -1)/ 3}$ must both
grow (or stay the same size), which is impossible taking into account
the first paragraph. Thus $(a,b,c) \in
S_a$.  But the triple $(a,b,c)$ was arbitrary, so every triple is in
$S_a$, and likewise in $S_b$ and $S_c$.\qed
\vs.1 \nd
Proof of Lemma 5.  First we show that for each horoball $a'-a$ must be
an integer.  For each dyadic integer $x$ let $\pi_m(x)$ be the
fractional part of $2^{-m}x$, and let $\sigma_m(x)$ be the integer
part.  The effect of adding $a'-a$ on the $m$th layer of the horoball
is to translate the locations of the tiles by $\pi_m(a'-a)$ and to
change the pattern of ``ancestor'' tiles by $\sigma_m(a'-a)$.  If
$a'-a$ is not an integer, then $\pi_m(a'-a)$ does not converge, so
different layers deep in the horoball get shifted by different
amounts, which contradicts uniform continuity, insofar as each piece of
each layer looks like a piece of every other layer.

Thus for each horoball, the difference $a'-a$ is an integer.  If these
differences are not bounded, we can pick a sequence of horoballs $H_m$
(with indices $a_m$, and possibly in different tilings) such that
$\pi_m(a'_m-a_m)$ does not converge.  Since with radius $m\log(2)$
every patch centered on the $m$th layer of $H_m$ is replicated in the
$M$th layer of $H_M$, for every $M>m$, this lack of convergence of
$\pi_m(a'_m-a_m)$ contradicts uniform continuity.\qed
\vs.1 \nd
Proof of Lemma 6.
Since the differences $(a'-a)$ take values in a finite set, the values
of $(a'-a)$ can be determined by knowing the first $N$ digits of
$(a,b,c)$.  But at some triangles, it happens that $a=b=c$\hs-.05  $\pmod{2^N}$
(since all allowable triples mod $2^N$ do occur, and since 3 is a unit
when working mod $2^N$).  At such triangles, we must have
$a'-a=b'-b=c'-c = (w'-w)/3$ by symmetry.  By Lemma 5 this common
difference must be a (rational) integer, so $w'-w \in 3\Z$, which completes
the proof of this lemma, and the theorem.\qed
\vs.1

A consequence of Lemma 6, together with the fact that two adjacent
horoballs determine the entire tiling, is:
\vs.1 \nd
{\bf Proposition 1}. There is a unique topological conjugacy from $\T_w$
to $\T_{w'}$ when $w-w' \in 3\Z$. Equivalently, there are no 
nontrivial automorphisms of $\T_w$.

\vs.2
Although the tilings in $\T(\bar \tau)$ -- in particular those in any
$\T_w$ -- cannot have a cofinite symmetry group, those in any $\T_w$ in
fact do have a nontrivial symmetry group, as we see next.

Let $T_0\in\T_w-\T_{deg}$
and consider the set $h^{-1}[h(T_0)]$ of nondegenerate tilings having the
same associated horoball packing as $T_0$.  We shall describe the
symmetry group of $h^{-1}[h(T_0)],$ i.e., the subgroup of PSL$(2,\R)$ 
consisting of those elements which fix every tiling in $h^{-1}[h(T_0)]$.  
Up to conjugacy, this group is independent of $w$ and $T_0$.

The symmetry group of the horoball packing $h(T_0)$ is conjugate in
PSL$(2,\R)$ to PSL$(2,\Z)$; choosing a particular conjugacy is the
same as choosing a triangle and a distinguished vertex in $h(T_0)$,
i.e., a state in the trivalent tree.  Fix a conjugacy and identify the 
symmetries of $h^{-1}[h(T_0)]$ with the elements of PSL$(2,\Z)$ which 
fix every pair of dyadic numbers. We recall some terminology from Lemma 2
concerning particular elements of PSL$(2,\Z)$. 
\vs.1 \nd
{\bf Theorem 6}.  The elements $L^2$ and $R^2\equiv C^2LC^2L$ freely
generate a subgroup $\F$ of PSL$(2,\Z)$ of index $6.$  The symmetry 
group of $h^{-1}[h(T_0)]$ is the kernel of the abelianization 
$\langle L^2,R^2\rangle\mapsto\Z\oplus\Z.$
\vs.1 \nd
Proof. The actions of $L^2$ and $R^2$ on pairs of dyadics are given by
$$
R^2:(a,b)\mapsto(a+1,b-2)\quad{\rm and}\quad L^2:(a,b)\mapsto(a+2,b-1). \cheqno\newprf
$$
One readily checks that the operations
 $L^2,R^2,L^{-2},R^{-2},R^2L^{-2}$
and $L^2R^{-2}$ are precisely the ones which move from a vertex in the 
trivalent tree to a vertex two edges away and induce maps of the form 
$(a,b)\mapsto(a+k,b+\ell)$, $k,\ell\in\Z$ on pairs of dyadic numbers.  
It follows that the index of $\F$ in PSL$(2,\Z)$ is $6.$

Freeness follows from the fact that distance from the starting point 
does not decrease as we follow some sequence of the basic operations 
$L^{\pm2},R^{\pm2}$ unless one of the operations is followed
immediately by its inverse.

Since the vectors $(1,-2)$ and $(2,-1)$ are linearly independent, 
the symmetry group of $h^{-1}[h(T_0)]$ consists of those words in 
$L^{\pm2},R^{\pm2}$ for which the sums of the powers of $R$ and $L$ 
are both zero, i.e., the kernel of the abelianization of $\F.$\qed
\vs.2
We now note that our use of the densest packing by horoballs,
Figure 6 (or, using the Poincar\'e disk, Figure 8), was not critical
in the above method. An infinite family of generalizations can be made
from other such horoball packings, as we now argue.

To generalize our ``triangular'' tilings, we consider tilings of 
$\H^2$ constructed as follows.
First pack $\H^2$ by horoballs such that five horoballs meet along
regular ``pentagons'' (rather than triangles), as in Figure 9. 
The symmetry group of such a packing, the Hecke group $G_5$, is 
a cofinite subgroup of PSL$(2,\R)$ generated by $z \to -1/z$ and 
$z \to z+ \lambda$, where $\lambda = (1 + \sqrt{5})/2$ is the golden mean. 

We
tile each horoball with (differently) modified binary tiles, where now
we need an appropriate width so we can
arrange that the prongs sticking out of such a tile each fill up a
fifth of a pentagon; see Figure 10. Relative to such a pentagon, the tiling of a
horoball with modified binary tiles is associated to a dyadic integer.
The first digit tells whether we are on the left or right prong of the
tile, the next digit tells whether that tile emerges from the left or
right prong of its parent, and so on.

Pick a rational integer $k$, once and for all. 
Let the dyadic integers around a pentagon, counting clockwise,
be $a$, $b$, $c$, $d$ and $e$. The
``$a$'' and ``$b$'' horoballs also meet at another pentagon,
and we assume that the five dyadic integers representing these
horoballs, counting {\it counterclockwise}, are $a+1$, $b -1$,
$c -k$, $d$ and $e+k$. This rule for relating patterns around adjacent
pentagons is a generalization of the ``fixed sum'' rule for triangular
tilings.  Notice also that the sum around the pentagons {\it is} fixed, 
and we let  $\T_{k,w}$ be
the closure of the class of tilings that follow the ``$k$-rule'' and
for which the sum of the five
dyadics at each pentagon add up to $w$. Again we denote by $\T_{deg}$
the degenerate tilings.
\vs.1 \nd
{\bf Theorem 7}. $\T_{k,w}$ is uniquely ergodic under the action
of $\G^2=\,$PSL$(2,\R)$.
\vs.1 \nd
Proof. This proof is nearly identical to the proof of Theorem 4. As
before, we approximate tilings in $\T_{k,w}$ by packings in which 
each horoball has only $N$ layers, and show that such packings have 
cofinite symmetry groups.  

As before let $\Theta_N(T)$ be the packing obtained from $T$ by removing all
but the $N$ 
horocyclic rows of tiles closest to the
boundary of any horoball in $h(T)$, and let $P_N$ be the range of $\Theta_N(T)$. 
In place of Lemma 2, we need to prove:  
\vs.1 \nd
{\bf Lemma 7}. For any $T\, T'  \in \T_{k,w} -\T_{deg}$, $\Theta_N(T)$ has
cofinite symmetry group, $\Theta_N(T)$ is in the orbit of 
$\Theta_N(T')$ and $\P_N$ is equal to this orbit union
the empty packing.
\vs.1 \nd
Proof of Lemma 7. A ``pentagonal'' horoball packing of $\H^2$
corresponds to an infinite 5-valent tree, with the pentagons of the
packing corresponding to the vertices of the tree.  We can navigate
around the tree with two fundamental operations.  If a ``state'' is a
vertex together with a choice of one of the five edges leading out
from that vertex, then the two operations are
$$ P = \hbox{Rotate counterclockwise by 72 degrees} \cheqno\penta $$
$$ L = \hbox{Go forwards to the next vertex and bear hard left.}\cheqno\pentb $$
Together these generate $G_5$. 
In terms of PSL$(2,\R)$, $P$ is the elliptic element $z \to 1/(\lambda-z)$
while $L$ is
the parabolic element $z \to z+ \lambda$.

We list the horoballs around a vertex in counterclockwise order, starting
with the one to the right of the chosen edge. We need only list the
first four of the five horoballs, since if their dyadic integers 
are $a$, $b$, $c$, and $d$, then the last one must be $e=w -a-b-c-d$.

The actions of $L$ and $P$ are easy to compute, namely:
$$L : (a,b,c,d) \to (a+1,e+k,d,c-k) \cheqno\pentc $$
$$P : (a,b,c,d) \to (e, a, b, c). \cheqno\pentd $$
It is not hard to check the following elements:
$$ L^2 : (a,b,c,d) \to (a+2, b-1+k, c-k, d-k) \cheqno\pente $$
$$ P L^2 P^4 : (a,b,c,d) \to (a+k-1, b+2, c-1+k, d-k) \cheqno\pentf $$
$$ P^2 L^2 P^3: (a,b,c,d) \to (a-k, b+k-1, c+2, d+k-1) \cheqno\pentg $$
$$ P^3 L^2 P^2: (a,b,c,d) \to (a-k, b-k, c+k-1, d+2). \cheqno\penth $$
Thus the possible values of $(a,b,c,d)$ differ by (among others) the elements
of the sub-lattice of $\Z^4$ generated by $(2, k-1, -k, -k)$, $(k-1,2,k-1,-k)$,
\vs0 \nd$(-k, k-1,2,k-1)$ and $(-k, -k, k-1, 2)$.  

Since 
$$\det \pmatrix{2 & k-1 & -k & -k \cr k-1 & 2 & k-1 & -k \cr
-k & k-1 & 2 & k-1 \cr -k & -k & k-1 & 2} = 5[(k-2)(k-1)k(k+1)+1] \cheqno\penti $$
is odd, 
$G_5$ acts transitively on the space of 
quadruples $(a,b,c,d)$\hs-.05  $\pmod{2^N}$, and the subgroup that preserves a given
quadruple  (and hence the first $N$ rows of each horoball in the tiling), 
is an index $2^{4N}$ subgroup
of $G_5$, and hence is a cofinite subgroup of PSL$(2,R)$. The remainder of
the lemma, and the theorem, follow as in triangular case.\qed
\vs.1
Next we consider the question of conjugacy for these systems.
Again, we just modify the argument that worked for triangle tilings.
\vs.1
\noindent {\bf Theorem 8}. $\T_{k,w}$ and $\T_{k,w'}$ are topologically 
conjugate if and only if $w-w' \in 5\Z$. $\T_w$ and $\T_{w'}$ are almost
conjugate for any $w,w'$.
\vs.1 \nd
Proof. The proof is essentially the same as the proof of Theorem 5, in
particular that of almost conjugacy, which we do not discuss further. If 
$w-w' \in 5\Z$, the conjugacy is simply adding $(w-w')/5$ to each 
dyadic index.  The converse follows from the analogues of Lemmas 3--6.
The proofs of Lemmas 3, 5, and 6 carry over almost word-for-word.  
Lemma 4 was algebraic, and used specific properties of the fixed-sum rule
for triangles.  In its place we have the following two lemmas that are
specific to pentagonal horoball packings and the Hecke group $G_5$.

\vs.1 \nd 
{\bf Lemma 8}. The only elements of $\Q[\sqrt{2}]$ 
that appear as elements of matrices in
$G_5$ are $-1$, $0$ and $1$. 
\vs.1 \nd
Proof. Viewed as matrices, the pentagonal Hecke group is generated by
$\pmatrix{1 & \lambda \cr 0 & 1}$ and $\pmatrix{0 & 1 \cr -1 & 0}$.  The
matrix elements are manifestly elements of $\Z[\lambda]$, and 
\vs0 \nd $\Z[\lambda]\cap \Q[\sqrt{2}]=\Z$. We will show that
the only integers that actually appear as matrix elements are 0 and
$\pm 1$.

Rosen [Rose] showed that an element of $\Z[\lambda]$ is an element of a
matrix in the group if and only if it is (up to sign) the denominator
of a finite approximant of the continued fraction 

$$r_0 \lambda + {\epsilon_1\over r_1 \lambda + {\epsilon_2\over r_2 \lambda +
\cdots}}, \cheqno\rosa $$

\nd where $\epsilon_i = \pm 1$ and each $r_i$ is a positive integer
(except possibly $r_0$, which may be zero). The continued fraction
expansion of a real number is not unique, but can always be expressed
in a unique ``reduced form'', one of whose requirements is that if
$r_n=1$, then $\epsilon_{n+1} r_{n+1} \ne -1$. The denominators $Q_n$
of the successive approximants to the (possibly infinite) continued
fraction satisfy the recursion:

$$ Q_n = r_n \lambda Q_{n-1} + \epsilon_n Q_{n-2}, \cheqno\rosb $$

\nd and we may take $Q_{-1}=0$ and $Q_0=1$.  Writing $Q_n = a_n \lambda
+ b_n$,
the recursion becomes:

$$ a_n = r_n (a_{n-1}+b_{n-1}) + \epsilon_n a_{n-2} \cheqno\rosc $$
$$ b_n = r_n a_{n-1} + \epsilon_n b_{n-2} \cheqno\rosd $$

\nd We claim that the coefficients satisfy three properties:
$\alpha$) $a_n \ge a_{n-1}$, $\beta$) $b_n \ge 0$ and $\gamma$) $a_n \ge b_{n-1}$.  These
are easily checked for $n=1, 2, 3$.  We prove these hold for all $n$
by induction.  Suppose they hold for $n$ up to $k$.
We have $a_{k+1} = r_{k+1}(a_k + b_k) +
\epsilon_{k+1}a_{k-1}$. If
$\epsilon_{k+1}=1$ this is manifestly at least $a_k$.  If $r_{k+1}>1$
and $\epsilon_{k+1}=-1$ then $a_{k+1} - a_k = r_{k+1} b_k +
(r_{k+1}-1) a_k - a_{k-1}$, which is non-negative since $a_k - a_{k-1}
\ge 0$.  Finally, if $r_{k+1}=1$ and $\epsilon_{k+1}=-1$, then $r_k
\ge 2$, and by property $\gamma$ we have
$$ a_{k+1} = a_k + b_k - a_{k-1}= a_k + (r_k-1)a_{k-1} +
\epsilon_k b_{k-2} \ge a_k + a_{k-1} - b_{k-2} \ge a_k, \cheqno\rose $$

\nd which is the needed induction for $\alpha$.
Next, $b_{k+1}  = r_{k+1} a_k + \epsilon_{k+1} b_{k-1}\ge a_k -
b_{k-1} \ge 0$, which is the needed induction for $\beta$.
For $\gamma$ we note $a_{k+1} \ge a_k + b_k - a_{k-1} \ge b_k$, which completes the induction.
Finally, since the sequence $a_k$ is nondecreasing and since
$a_1$ is positive $a_k$ is never zero and $Q_k$ is never rational
for $k\ge 1$.\qed 
\vs.1 \nd
{\bf Lemma 9}.
A conjugacy between $\T_{k,w}$ and $\T_{k,w'}$
must preserve the locations of horoballs exactly.
\vs.1 \nd 
Proof. We have already shown that the points on the sphere at
infinity where the horoballs touch are not changed, that their
radii change by a bounded amount, and that there are only a
finite number of possible values for that change in radius. This
implies that knowing the first $N$ digits of all five indices at
a pentagon will determine which horoballs grow and shrink, and by
how much. Note also that the deep interiors of all horoballs are
identical, so the change in radius is the same for all
horoballs, modulo $\log(2)$.

Now consider a tiling whose associated horoball packing is as follows:
One horoball is the set $\{ x+iy\,|\, y \ge 1\}$, and the others are its
images under the Hecke group $G_5$.
For each matrix $\pmatrix{\alpha & \beta \cr \gamma & \delta}$ in the
group there is a horoball tangent to the $x$-axis at
$\alpha/\gamma$ with Euclidean diameter $1/\gamma^2$.

Note that the dyadic indices of the tiling modulo $2^N$ are unchanged
by the transformation $z \to z+ 2^{N+1}\lambda$, so that $\phi$ of this
tiling corresponds to a packing that is invariant under $z \to z+
2^{N+1} \lambda$.  However, if the packing contains the horoball $\{x+iy\,|\,
y \ge \nu\}$, then it is invariant only under addition of multiples of
$\nu \lambda$.  Thus $\nu$ must divide $2^{N+1}$, and in particular must
be rational.

Now consider a horoball that meets the horoball at infinity in the new
packing (but did not meet the horoball at infinity in the original
packing).  Since its (hyperbolic) radius has changed by $-\log(\nu)$ \hs-.15 
$\pmod{\log(2)}$, and since $\nu$ is rational, its Euclidean
radius must have been a power of 2 times the square of a rational to
begin with.  However, this implies that there is an element of the
Hecke group with $\gamma$ of the form of a product of a rational and
a power of $\sqrt{2}$, and not equal to 0 or $\pm 1$. By Lemma 8 no such
element exists.\qed
\vs.1
Next we consider some differences between optimization problems with
other variations on our basic tile.  So far we have considered
packings of $\H^2$ by horoballs which meet either in ``triangles''
(the densest packing of horoballs) or in ``pentagons'', and modified
our basic tile to have 2 prongs, each of which is either a third of a
triangle (Figure 7) or a fifth of a pentagon (Figure 10). One can
easily allow horoball packings defined by other regular $n$-gons, and
also consider tiles to have more prongs, one for each neighboring
$n$-gon; for $n=3$ and $m=3$ see Figure 11.

Let $\T({n,m})$ denote the space of all tilings by $m$-pronged tiles,
each prong congruent to one of the $n$ isosceles triangles dividing an
$n$-gon as defined above.  Let $X$ be a closed, PSL$(2,\R)$-invariant
subset of $\T({n,m})$ and let $Y$ be a closed, PSL$(2,\R)$-invariant
subset of $\T({n',m'})$.  We want to show that, under various
assumptions, $X$ and $Y$ cannot be topologically conjugate.
\vs.1
\nd {\bf Lemma 10}. If $\phi: X \to Y$ is a topological conjugacy, 
$\phi$ maps degenerate tilings to degenerate tilings and
nondegenerate tilings to nondegenerate tilings. Moreover, for
nondegenerate tilings the points of tangency of the horoballs at the
sphere at infinity are not changed by $\phi$.
\vs.1 \nd
Proof. The proof is essentially that of Lemma 3. The fact that we were 
dealing with fixed-sum tilings (with $n=n'=3$ and $m=m'=2$) was never used. 
\vs.1 \nd 
{\bf Theorem 9}. Let $n=3$ and $n'=5$. Suppose $X\subseteq \T({n,m})$
contains nondegenerate tilings, is invariant under $\G^n$ and closed,
and suppose $Y\subseteq \T({n',m'})$ is invariant under $\G^n$ and
closed. Then the actions of $\G^n$ on $X$ and $Y$ are not
topologically conjugate.
\vs.1 \nd
Proof. For a nondegenerate tiling in $\T({3,m})$, the ``cusp point set'',
the set of points at infinity that meet horoballs, is
conjugate, by some fixed element of PSL$(2,\R)$ to $\Q \cup \{\infty\}$.
The cusp point set for the corresponding (nondegenerate) tiling 
in $\T({5,m'})$ is conjugate to $\Q[\lambda] \cup \infty$, 
where $\lambda$ is the golden mean. This contradicts Lemma 10.\qed
\vs.1 \nd
{\bf Theorem 10}. Suppose $X\subseteq \T({n,m})$ is invariant under
$\G^n$ and closed, and $Y\subseteq \T({n',m'})$ is invariant under
$\G^n$ and closed. If $m\ne m'$ then the actions of $\G^n$ on $X$ and
$Y$ are not topologically conjugate.
\vs.1 \nd
Proof: Assume without loss of generality that $m<m'$.
$X$ necessarily 
contains all the degenerate tilings, and in particular contains a tiling 
invariant under the map $z \to mz$.  By Lemma 10, any conjugacy would 
have to take this to a degenerate tiling in $\T({n',m'})$ that is also 
invariant under $z \to mz$.  However, under this symmetry, points along 
the $y$-axis are only moved a distance $\ln(m)$, while any symmetry of a 
degenerate tiling in
$\T({n',m'})$ must move points at least a distance $\ln(m')$.\qed
\vs.1
Note that Theorem 9 requires that $X$ contain nondegenerate tilings;
if $X$ consists only of degenerate tilings, then $n$ is irrelevant.
Theorem 10, however, only depends on the existence of degenerate
tilings. These arise automatically since they are in the orbit closure of
every nondegenerate tiling. 
Generalizing Theorem 9 would require knowing more than we do
about cusp point sets for Hecke groups.
\vs.2 \nd
{\bf V. Isomorphism and Uniqueness in Problems of Optimally Dense}
\vs0 \nd {\bf Packings}
\vs.1 
Optimally dense packings, especially tilings of $\E^2$ and $\E^3$,
have been important for many years in classifying certain geometric
properties of patterns; we refer here to the classification through
symmetry groups, the so-called crystallographic symmetries.

At heart the formalism consists of treating the structures of
interest as subsets of $\S^n$, with the action on them of $\G^n$ --
that is, one uses the structure of dynamical systems.  Consider for
instance two ball packings in $\E^3$, $P^{fcc}$, the face centered
cubic, and $P^{hcp}$, the hexagonal close packed. When we choose to
distinguish $P^{fcc}$ from $P^{hcp}$ on symmetry grounds what we are
saying is that the subgroup of $\G^3$ (the connected Euclidean group)
which acts trivially on every element of one orbit is different from
the symmetry group of the elements of the other orbit. This implies
the systems are not conjugate: there is no bijection, between the two
orbits under $\G^3$, which intertwines the action of $\G^3$. 
So in this simple situation we see that conjugacy can
detect differences in symmetry.

The case of the kite \& dart tilings (Figure 3) is instructive. As we noted in
section III, it is natural to want to think of all these tilings as
equivalent. This is true even though the tilings can actually have
different symmetry groups; for instance there are two noncongruent kite
\& dart tilings with a point of 5-fold rotational symmetry, which the
other kite \& dart tilings do not have [Gard]. Furthermore, the symmetry of
the two special tilings actually do not play an essential role, for
two reasons. First, the 5-fold rotational symmetry appears in
regions of arbitrarily large size 
in every tiling, and this could replace the exact
symmetry of the special tilings. Furthermore the Penrose
tilings have a statistical form of 10-fold rotational symmetry,
which is expressed by the 10-fold rotational symmetry of all the
translation invariant measures on the space of Penrose tilings [Radi].
We also note an analogy between the different ``symmetry'' of elements of
tilings with different fixed-sums, as evidenced by nonconjugacy, and
the different (5-fold rotational) symmetry that appears among kite \&
dart tilings.

So we are led to relax the strict form of equivalence whereby two
packings are equivalent if they are in the same orbit under $\G^n$.
 From the example of the kite \& dart tilings one might be led to
replace this by having optimally dense packings equivalent if they are
generic for the same measure. But from the various examples
of section IV we will go one step further.
\vs.1 \nd
{\bf Definition 4}. We say that two optimally dense packings
$p,p'$ are ``weakly equivalent'' if the optimal measures $\mu_p$,
$\mu_{p'}$ for which they are generic have the following property: the
set $M(p)$ of all optimal measures the support of which intersects
the support of $\mu_p$ coincides with the set $M(p')$. An optimal
density problem will be said to have a ``unique solution'' if there is
only one weak equivalence class of optimally dense packings.
\vs.1
As we saw in section III, there are simple examples of optimally dense
packing problems, in particular that for disks of fixed radius in
$\E^2$, for which the solution is unique in the sense of consisting of
a single orbit of $\G^n$, or, put another way, in the sense that the
quotient $P^o_\B/\G^n$ consists of a single point. Because of the
aperiodicity of the kite \& dart tilings, and the modified binary
tilings, we have been led to divide $P^o_\B$ by a cruder
equivalence relation.

This paper is a continuation of a long tradition of classifying a
pattern through the dynamical system associated with it by the action
of the isometry group of the ambient space of the pattern. A common
step taken when following the dynamics approach is to settle on a form
of conjugacy, typically either measurable conjugacy or topological
conjugacy -- and in effect declare two optimization (or tiling)
problems equivalent if their dynamical systems are conjugate in the
chosen sense. Prominent in this vein is the analysis by Connes,
Putnam, Kellendonk et al. noted above, in which invariants of
aperiodic tilings are sought through operator algebras associated with
their dynamical systems.
\vs.1 \nd
{\bf Definition 5}. We declare two optimal density problems,
associated with finite sets $\B$ and $\B'$ of bodies in some fixed
$\S^n$, to be ``equivalent'' if there is a topological conjugacy
between their dynamical systems, $(\bar {P}^{o}_\B,\G^n)$ and
$(\bar {P}^{o}_\B,\G^n)$, where $\bar {P}^{o}_\B$ denotes the
closure in ${P_\B}$ of ${P^o_\B}$.
\vs.1
In terms of this notion of equivalence the proofs of Theorems 9 and 10
show how the optimization problems, for different variations of our
modified binary tile, can be distinguished by geometric features.

The operator algebra approach noted above is a powerful way to obtain
the desired invariants for topological conjugacy. Associated to a
dynamical system $(X, G)$ is the crossed-product $C^*$-algebra
$C(X)\times_\alpha G$, where $C(X)$ is the $C^*$-algebra of continuous
complex-valued functions on $X$ and $\alpha$ is the action of $G$ on
$C(X)$. (The crossed-product is the completion of the algebraic tensor
product in a certain norm; for this and other terms in operator
algebras we refer to [Blac].) The $K$-theoretic invariants of the
crossed-product algebra are topological conjugacy invariants for the
dynamical system.

A common way to compute invariants is to associate an AF algebra with
the dynamical system, ``large'' in an algebra Morita equivalent to the
crossed-product of interest. (Morita equivalence preserves the
$K$-theoretic invariants.) We do not see how to do that
here. Alternatively one could try to compute $K_0$ using tools such as the Pimsner-Voiculescu
6-term exact sequence; this has
been a practical route at least when the dynamical group is $\R$ or
$\R^2$, but this seems to be harder for groups such as PSL$(2,\R)$.

In short, the operator algebra methods used to produce dynamical invariants for
aperiodic systems in Euclidean space seem to need extension for this
more general aperiodic setting.

\vs2 \nd
{\bf VI. Summary} 
\vs.1
We have considered in a unified fashion, for Euclidean and hyperbolic
spaces, the class of problems of determining the optimally dense
packings by congruent copies of bodies from some fixed finite set
$\B$. 

 From a qualitative point of view the main solved examples in the
literature are as follows. For balls of fixed radius in $\E^2$ or
fixed tight radius in $\H^2$ or $\H^4$ there is a unique periodic
solution [BoR1]. For balls of fixed radius in $\E^3$ the problem is widely
believed degenerate, with nonunique solution, but including periodic
special cases. For the general problem for balls of fixed radius the
problem is better understood in $\H^n$ than in $\E^n$, at least in the
qualitative sense we are emphasizing; for tight radius the problem has
the obvious unique periodic solution, and for most radii (all but
countably many) the problem in $\H^n$ is aperiodic, that is, {\it
none} of the optimally dense packings have symmetry group with
fundamental domain of finite volume [BoR1]. (For non-tight radii none
of the optimal packings have been published.) For optimal density
problems of polyhedra there is again the phenomenon of aperiodicity,
exemplified by the kite \& dart tilings in $\E^2$. Aperiodicity was
discovered in the Euclidean context by Berger in 1966 [Berg], and has
been of growing interest, from many perspectives (though mainly still
in Euclidean spaces) ever since [Radi].

We emphasize how short the above list is, how few optimum density
problems have been solved. As part of this small but distinguished
class, the aperiodic ones must play a significant role in any general
understanding of optimum density packings.

To prove existence of optimally dense packings it has been found
useful to work in an ergodic theoretic formalism; in a sense, the
solutions of such problems are naturally organized through invariant
probability measures on a topological space of packings. Much of this
paper is an analysis within that formalism of certain new examples:
the densest packings of the hyperbolic plane by congruent copies of a
certain polygon, the modified binary $\bar \tau$, and others in a
2-parameter family of variations. These are the first aperiodic
optimal density problems in hyperbolic space for which explicit solutions
have been determined, and they exhibit features unknown from
examples set in Euclidean spaces -- specifically, the small, compact,
invariant subset $\T_{deg}$ appearing in the supports of all invariant
measures on the solution space.

Optimal density problems have a long tradition but have never been
treated as have other classes of optimization problems -- for instance
by analysis of conditions for existence and uniqueness of
solutions. In this paper we consider the question of uniqueness,
guided by symmetry properties of solutions. We analyze such problems
as dynamical systems, with the group of isometries of the ambient
space ($\E^n$ or $\H^n$) acting on the (compact, metrizable) space of
all possible packings. We are led to classify such systems up to
topological conjugacy, and use geometric features as invariants.

More specifically, Theorems 1 and 2 led us to partition the class of
optimal density problems with unique solution into two classes: the
periodic and the aperiodic. The study of {\it periodic} optimization
problems in $\E^n$, for which the solutions have cofinite symmetry
group, led, many years ago, to classification of the discrete
subgroups of the isometry group of Euclidean space in low
dimensions. Aperiodic tilings have led to related work; among problems
set in $\E^3$, the study of the quaquaversal aperiodic tilings [CoRa]
led to new results on classification of certain (dense) subgroups of
$SO(3)$ [RaS1], [RaS2], [CoRS]. Similarly, the classification results
in Section IV of our optimization problems in $\H^2$ naturally led to
{\it noncofinite} subgroups of PSL$(2\R)$, such as the symmetry groups
of fixed-sum tilings (Theorem 6), as well as questions about Hecke groups.

In summary, classification of aperiodic optimization problems amounts
to a study of the ``symmetries'' of the packing solutions in senses
related to, but different from, the manner appropriate for the well
studied periodic structures. The mathematics that is generated by such
analysis, a mixture of dynamics, operator algebras and Lie groups, is 
perhaps the main significance of the study of optimal density problems.
\vs.2 \nd
{\bf Acknowledgments}. We are grateful for useful discussions with Alan Reid, 
in particular for pointing us to reference [Rose].
\vfill \eject

\centerline{{\bf References}}
\vs.2 \nd

\item{[AdMa]} R. Adler and B. Marcus, Topological entropy and equivalence
of dynamical systems, {\it Mem. Amer. Math. Soc.} 219(1979).

\item{[Bear]} A. Beardon, {\it The Geometry of Discrete Groups},
Springer-Verlag, New York, 1983.

\item{[Berg]}  R. Berger, The undecidability of the domino problem, {\it 
Mem. Amer. Math. Soc.} 66(1966).

\item{[Blac]} B. Blackadar, {\it $K$-Theory for Operator Algebras}, 2nd ed,
Cambridge University Press, Cambridge, 1998.

\item{[Boro]} K. B\"or\"oczky, Packing of spheres in spaces of constant
curvature, {\it Acta Math. Acad.  Sci. Hung.} 32(1978), 243--261.

\item{[Bowe]} L. Bowen, On the existence of completely saturated packings and
completely reduced coverings, {\it Geometria Dedicata} 98(2003), 211--226.

\item{[BoR1]} L. Bowen and C. Radin, Densest packing of equal spheres in
hyperbolic space, {\it Discrete Comput. Geom.} 29(2003), 23--39.

\item{[BoR2]} L. Bowen and C. Radin, Optimally dense packings of hyperbolic space, 
{\it Geometriae Dedicata} 104(2004), 37--59.

\item{[Conn]} A. Connes, {\it Noncommutative Geometry}, Academic Press, San Diego, 1994.

\item{[CoRa]} J.H. Conway and C. Radin, Quaquaversal tilings and rotations, 
{\it Inventiones math.} 132(1998), 179--188.

\item{[CoRS]} J.H. Conway, C. Radin and L. Sadun, Relations in SO(3) supported by 
geodetic angles, {\it Discrete Comput. Geom.} 23(2000), 453--463.

\item{[CoGS]} J.H. Conway, C. Goodman-Strauss and N. Sloane, 
Recent progress in sphere packing,  pp. 37--76, in
{\it Current Developments in Mathematics}, Int. Press, Cambridge, 1999.
Somerville, MA, 1999.

\item{[Feje]} L. Fejes T\'oth, {\it Regular Figures}, Macmillan, New York,
1964.

\item{[FeKu]} G. Fejes T\'oth and W. Kuperberg, Packing and covering
with convex sets, chapter 3.3, pp. 799--860, in Vol B of {\it
Handbook of Convex Geometry}, ed.\ P. Gruber and J. Wills, North
Holland, Amsterdam, 1993.

\item{[Gard]} M. Gardner, Mathematical Games, in {\it Sci.\ Amer.}\ (January 1977) 
110--121.

\item{[Good]} C. Goodman-Strauss, A strongly aperiodic set of tiles in the 
hyperbolic plane, preprint, Univ. Arkansas, 2000, available from

\hs.5 http://comp.uark.edu/$\sim$cgstraus/papers/hyp.pdf

\item{[HF]} T. Hales, The Kepler conjecture, available from math.MG/9811078.

\item{[Hilb]} D. Hilbert, Mathematische Probleme, {\it Archiv der Mathematik und Physik}, 
(3) 1(1901), 44--63 and 213--237. Translated in {\it Bull. Amer. Math. Soc.} 8(1902), 437--479,
by Maby Winton Newson.

\item{[Hsia]} W.Y. Hsiang, On the sphere packing problem and the proof of
Kepler's conjecture, {\it Internat. J. Math.} 4(1993), 739--831.

\item {[Kato]} S. Katok, {\it Fuchsian Groups}, University of Chicago Press,
Chicago, 1992.
\vs1
\item{[KePu]} J. Kellendonk and I.F. Putnam, Tilings, C$^{\ast}$-algebras
and K-theory, in {\it Directions in Mathematical Quasicrystals}, ed. M. Baake
and R.V. Moody, CRM Monograph Series, Amer. Math. Soc., Providence, RI, 2000.

\item{[KuKu]} G. Kuperberg and W. Kuperberg, Double-lattice packings of convex 
bodies in the plane, {\it Discrete Comput. Geom.} 5(1990), 389--397.

\item{[Kupe]} G. Kuperberg, Notions of denseness, {\it Geom. Topol.}
4(2000) 277--292.

\item{[Laga]} J.C. Lagarias, Bounds for local density of sphere packings and
the Kepler conjecture, {\it Discrete Comput. Geom.} 27(2002) 165--193.

\item{[MaMo]} G.A. Margulis and S. Mozes, Aperiodic tilings of the
hyperbolic plane by convex polygons, {\it Israel J. Math.}
107(1998) 319--332.

\item{[MiRa]} J. Mi\c ekisz and C. Radin, Why solids are not really crystalline, 
{\it Phys. Rev.} B39(1989), 1950--1952. 

\item{[Miln]} J. Milnor, Hilbert's problem 18: On crystallographic groups,
fundamental domains, and on sphere packing, {\it Proc. Sympos. Pure
Math.} 28(1976), 491--506.

\item{[Moze]} S. Mozes, Aperiodic tilings, {\it Invent. Math.} 128(1997), 603--611. 

\item{[Nevo]} A. Nevo, Pointwise ergodic theorems for radial averages
on simple Lie groups I, {\it Duke Math. J.} 76(1994), 113--140.

\item{[NeSt]} A. Nevo and E. Stein, Analogs of Wiener's ergodic
theorems for semisimple groups I, {\it Annals of Math.} 145(1997),
565--595.

\item{[Penr]} R. Penrose, Pentaplexity - a class of non-periodic tilings of the
plane, {\it Eureka} 39(1978,) 16--22. (Reproduced in {\it Math. Intell.} 2(1979/80), 32--37.)

\item{[Radi]} C. Radin, {\it Miles of Tiles}, Student Mathematical Library,
Vol. 1,
Amer. Math. Soc., Providence, 1999.

\item{[RaS1]} C. Radin and L. Sadun, Subgroups of $SO(3)$ associated 
with tilings, {\it J. Algebra} 202(1998), 611--633.

\item{[RaS2]} C. Radin and L. Sadun, On 2-generator subgroups of SO(3), 
{\it Trans. Amer. Math. Soc.} 351(1999), 4469--4480.

\item{[RaWo]} C. Radin and M. Wolff, Space tilings and local isomorphism,
{\it Geometriae Dedicata} 42(1992), 355--360.

\item{[Roge]} C.A. Rogers, {\it Packing and Covering}, University Press,
Cambridge, 1964.

\item{[Rose]} D. Rosen, A class of continued fractions associated with certain
properly discontinuous groups, {\it Duke Math. J.} 21(1954),
549--563.

\item{[Walt]} P. Walters, {\it An Introduction to Ergodic Theory},
Springer-Verlag, New York, 1982.
\vfill \eject


%
%
\newdimen\FigSize       \FigSize=.9\hsize 
%
\newskip\abovefigskip   \newskip\belowfigskip
\gdef\epsfig#1;#2;{\par\vskip\abovefigskip\penalty -500
   {\everypar={}\epsfxsize=#1\nd \centerline{\epsfbox{#2}}}%
    \vskip\belowfigskip}%
%
\newskip\figtitleskip
\gdef\tepsfig#1;#2;#3{\par\vskip\abovefigskip\penalty -500
   {\everypar={}\epsfxsize=#1\nd
    \vbox
      {\centerline{\epsfbox{#2}}\vskip\figtitleskip
       \centerline{\figtitlefont#3}}}%
    \vskip\belowfigskip}%
%
\newcount\FigNr \global\FigNr=0
\gdef\nepsfig#1;#2;#3{\global\advance\FigNr by 1
   \tepsfig#1;#2;{Figure\space\the\FigNr.\space#3}}%
%
%
%
\gdef\ipsfig#1;#2;{
   \midinsert{\everypar={}\epsfxsize=#1\nd
              \centerline{\epsfbox{#2}}}
   \endinsert}%
%
\gdef\tipsfig#1;#2;#3{\midinsert
   {\everypar={}\epsfxsize=#1\nd
    \vbox{\centerline{\epsfbox{#2}}%
          \vskip\figtitleskip
          \centerline{\figtitlefont#3}}}\endinsert}%
%
\gdef\nipsfig#1;#2;#3{\global\advance\FigNr by1%
  \tipsfig#1;#2;{Figure\space\the\FigNr.\space#3}}%
\newread\epsffilein    
\newif\ifepsffileok    
\newif\ifepsfbbfound   
\newif\ifepsfverbose   
\newdimen\epsfxsize    
\newdimen\epsfysize    
\newdimen\epsftsize    
\newdimen\epsfrsize    
\newdimen\epsftmp      
\newdimen\pspoints     
\pspoints=1bp          
\epsfxsize=0pt         
\epsfysize=0pt         
\def\epsfbox#1{\global\def\epsfllx{72}\global\def\epsflly{72}%
   \global\def\epsfurx{540}\global\def\epsfury{720}%
   \def\lbracket{[}\def\testit{#1}\ifx\testit\lbracket
   \let\next=\epsfgetlitbb\else\let\next=\epsfnormal\fi\next{#1}}%
\def\epsfgetlitbb#1#2 #3 #4 #5]#6{\epsfgrab #2 #3 #4 #5 .\\%
   \epsfsetgraph{#6}}%
\def\epsfnormal#1{\epsfgetbb{#1}\epsfsetgraph{#1}}%
\def\epsfgetbb#1{%
%
%
\openin\epsffilein=#1
\ifeof\epsffilein\errmessage{I couldn't open #1, will ignore it}\else
%
%
   {\epsffileoktrue \chardef\other=12
    \def\do##1{\catcode`##1=\other}\dospecials \catcode`\ =10
    \loop
       \read\epsffilein to \epsffileline
       \ifeof\epsffilein\epsffileokfalse\else
%
%
          \expandafter\epsfaux\epsffileline:. \\%
       \fi
   \ifepsffileok\repeat
   \ifepsfbbfound\else
    \ifepsfverbose\message{No bounding box comment in #1; using
defaults}\fi\fi
   }\closein\epsffilein\fi}%
%
%
\def\epsfsetgraph#1{%
   \epsfrsize=\epsfury\pspoints
   \advance\epsfrsize by-\epsflly\pspoints
   \epsftsize=\epsfurx\pspoints
   \advance\epsftsize by-\epsfllx\pspoints
%
%
   \epsfxsize\epsfsize\epsftsize\epsfrsize
   \ifnum\epsfxsize=0 \ifnum\epsfysize=0
      \epsfxsize=\epsftsize \epsfysize=\epsfrsize
%
arithmetic!
%
     \else\epsftmp=\epsftsize \divide\epsftmp\epsfrsize
       \epsfxsize=\epsfysize \multiply\epsfxsize\epsftmp
       \multiply\epsftmp\epsfrsize \advance\epsftsize-\epsftmp
       \epsftmp=\epsfysize
       \loop \advance\epsftsize\epsftsize \divide\epsftmp 2
       \ifnum\epsftmp>0
          \ifnum\epsftsize<\epsfrsize\else
             \advance\epsftsize-\epsfrsize \advance\epsfxsize\epsftmp
\fi
       \repeat
     \fi
   \else\epsftmp=\epsfrsize \divide\epsftmp\epsftsize
     \epsfysize=\epsfxsize \multiply\epsfysize\epsftmp
     \multiply\epsftmp\epsftsize \advance\epsfrsize-\epsftmp
     \epsftmp=\epsfxsize
     \loop \advance\epsfrsize\epsfrsize \divide\epsftmp 2
     \ifnum\epsftmp>0
        \ifnum\epsfrsize<\epsftsize\else
           \advance\epsfrsize-\epsftsize \advance\epsfysize\epsftmp \fi
     \repeat
   \fi
%
%
   \ifepsfverbose\message{#1: width=\the\epsfxsize,
height=\the\epsfysize}\fi
   \epsftmp=10\epsfxsize \divide\epsftmp\pspoints
   \vbox to\epsfysize{\vfil\hbox to\epsfxsize{%
      \includegraphics{#1}%
      \hfil}}%
\epsfxsize=0pt\epsfysize=0pt}%
%
%
{\catcode`\%=12
\global\let\epsfpercent=
%
%
\long\def\epsfaux#1#2:#3\\{\ifx#1\epsfpercent
   \def\testit{#2}\ifx\testit\epsfbblit
      \epsfgrab #3 . . . \\%
      \epsffileokfalse
      \global\epsfbbfoundtrue
   \fi\else\ifx#1\par\else\epsffileokfalse\fi\fi}%
%
%
\def\epsfgrab #1 #2 #3 #4 #5\\{%
   \global\def\epsfllx{#1}\ifx\epsfllx\empty
      \epsfgrab #2 #3 #4 #5 .\\\else
   \global\def\epsflly{#2}%
   \global\def\epsfurx{#3}\global\def\epsfury{#4}\fi}%
%
%
\def\epsfsize#1#2{\epsfxsize}%
%
%

\epsfverbosetrue                        
\abovefigskip=\baselineskip             
\belowfigskip=\baselineskip             
\global\let\figtitlefont\bf             
\global\figtitleskip=.5\baselineskip    

\hbox{}\vs1
\vbox{\epsfig 1\hsize; 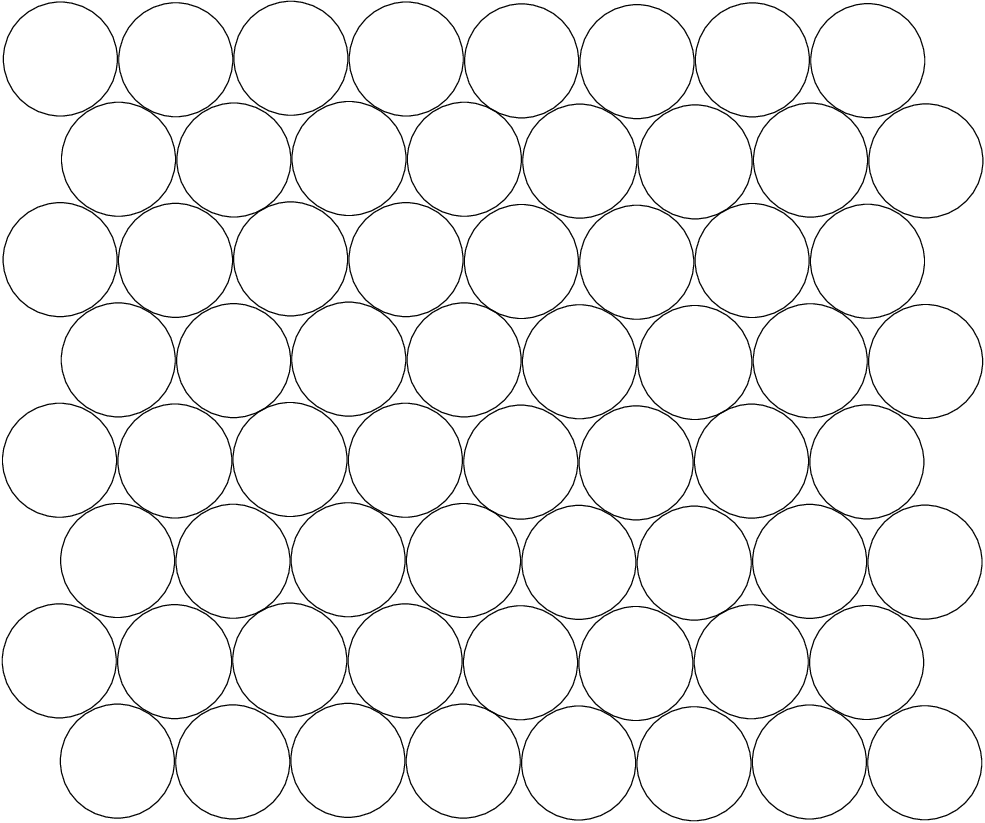 ;}
\vs.5
\centerline{Figure 1. The densest packing of equal disks
in the Euclidean plane}
\vfill \eject

\hbox{}\vs2 \nd
\epsfig 1\hsize; 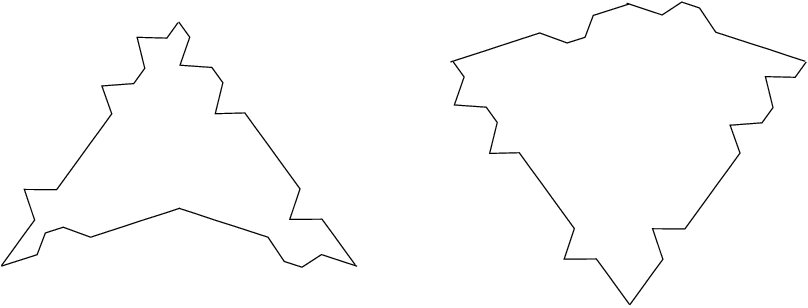;
\vs.3 \nd
\centerline{Figure 2. The Penrose kite \& dart tiles of the Euclidean plane}
\vfill \eject

\hbox{}\vs-.8
\epsfig 1.15\hsize; 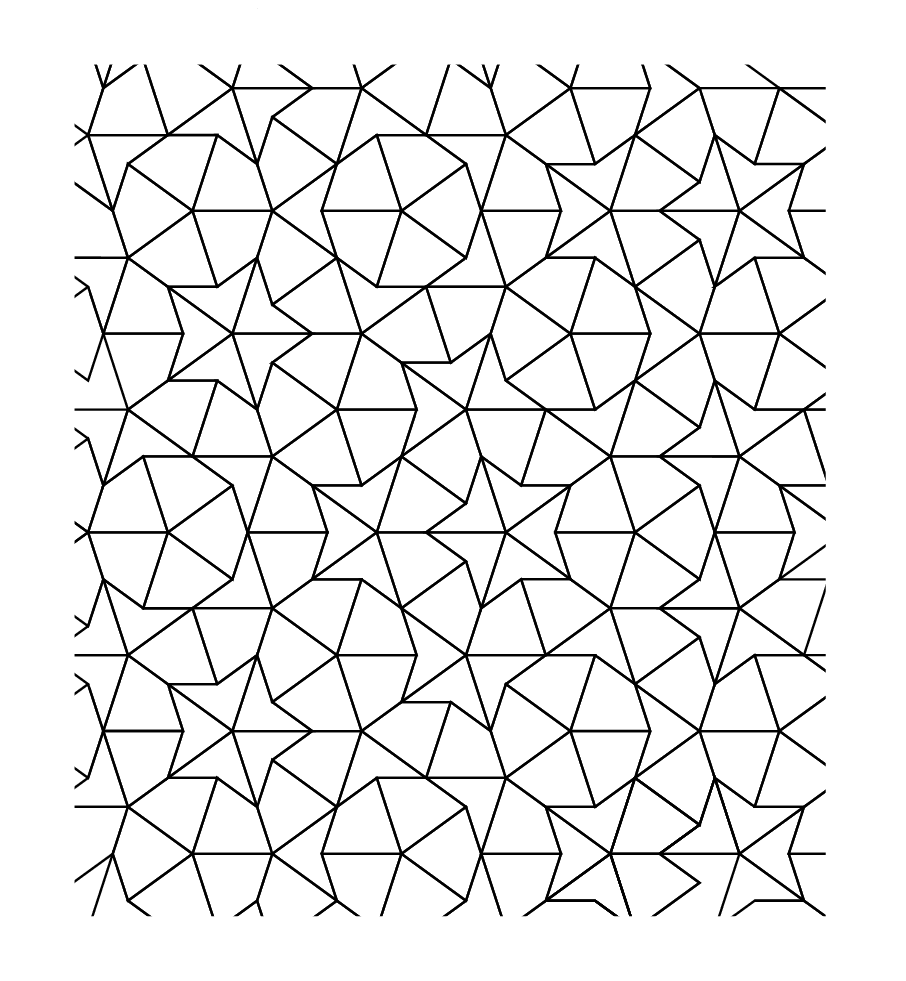 ;
\vs-.5
\centerline{Figure 3. A Penrose kite \& dart tiling of the Euclidean plane}
\vfill\eject

\hbox{}\vs1 \nd
\epsfig .8\hsize; 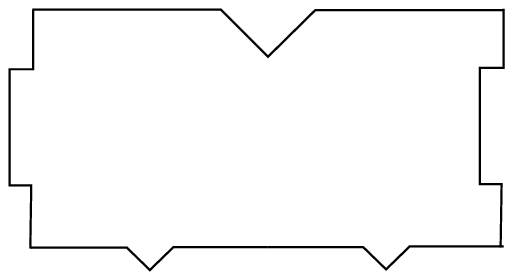;
\vs.3 \nd
\centerline{Figure 4. The simple binary tile in the upper half plane model of the hyperbolic plane}
\vfill\eject

\hbox{}\vs1
\vbox{\epsfig 1\hsize; 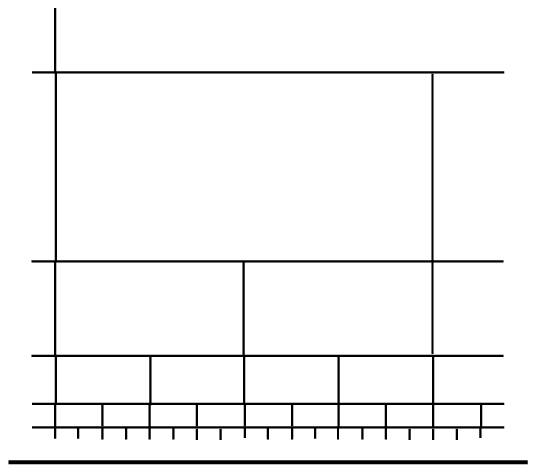 ;}
\vs.5
\centerline{Figure 5. The binary tiling of the upper half plane}
\vfill
\vfill \eject

\hbox{}
\vbox{\epsfig 1\hsize; 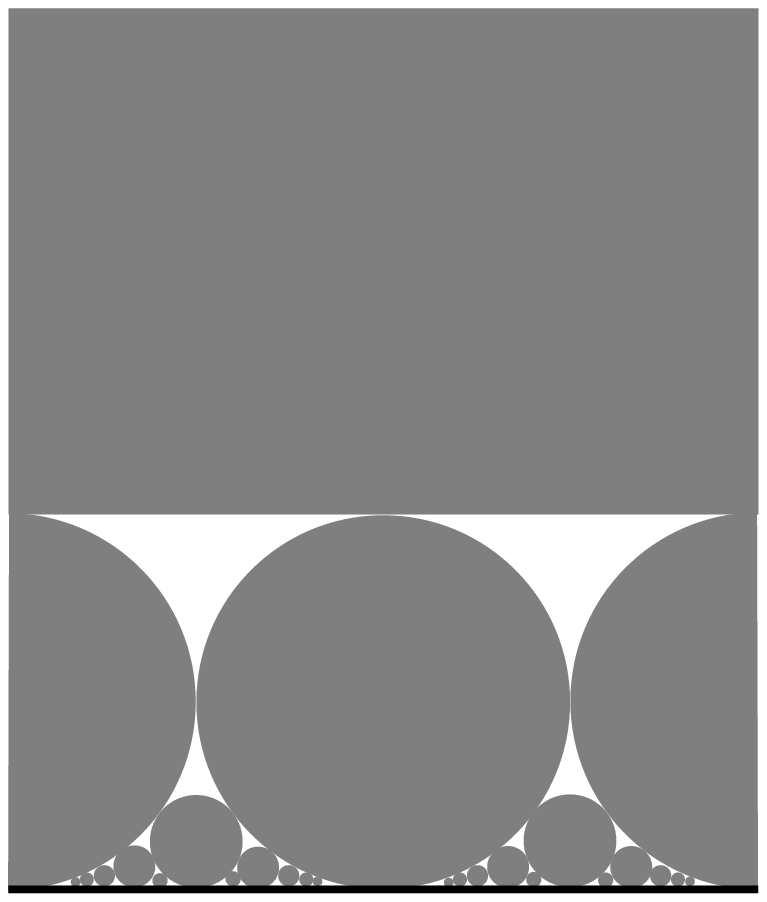;}
\vs.01
\centerline{Figure 6. A packing of the upper half plane by horoballs}
\vfill\eject

\hbox{}\vs1 \nd
\epsfig .8\hsize; 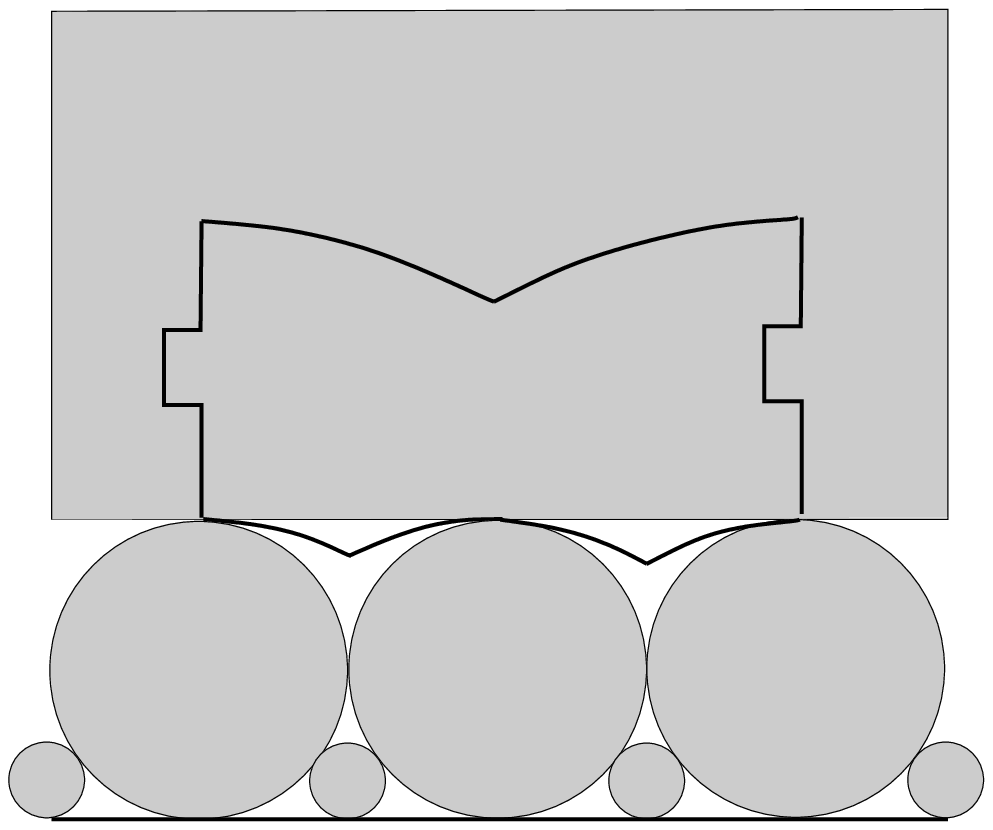;
\vs.3 \nd
\centerline{Figure 7. The modified binary tile in the upper half plane}
\vfill\eject

\hbox{}\vs1 \nd
\epsfig .8\hsize; 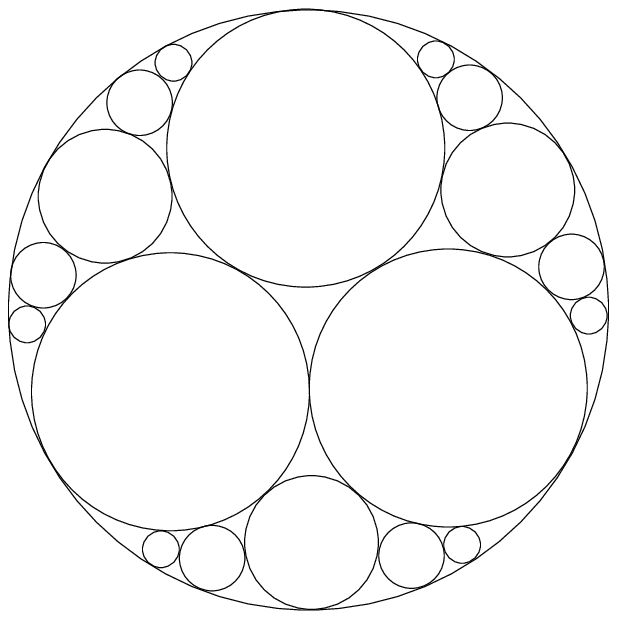;
\vs.3 \nd
\centerline{Figure 8. A packing of the Poincar\'e disk by horoballs}
\vfill\eject

\hbox{}\vs1 \nd
\epsfig .8\hsize; 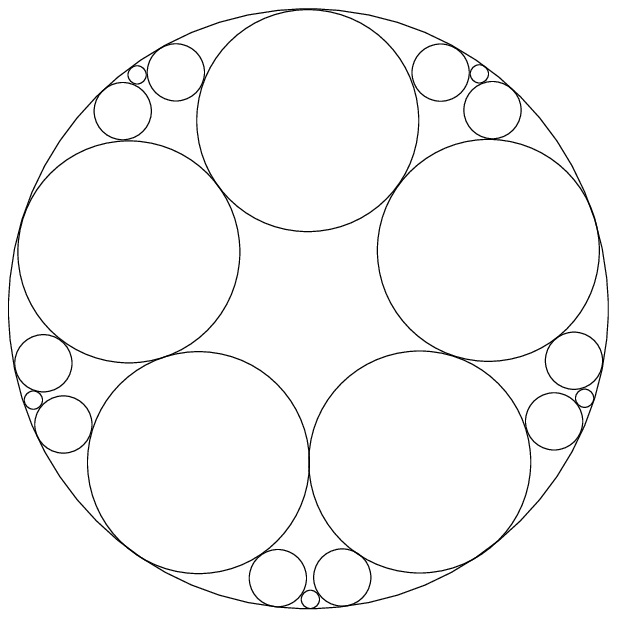;
\vs.3 \nd
\centerline{Figure 9. Another packing of the Poincar\'e disk by horoballs}
\vfill\eject

\hbox{}\vs1 \nd
\epsfig .8\hsize; 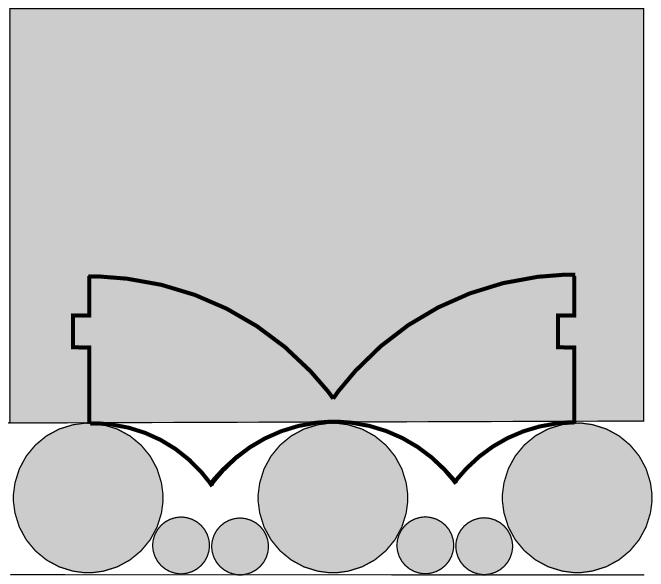;
\vs.3 \nd
\centerline{Figure 10. A ``pentagonal'' binary tile}
\vfill\eject

\hbox{}\vs1 \nd
\epsfig .8\hsize; 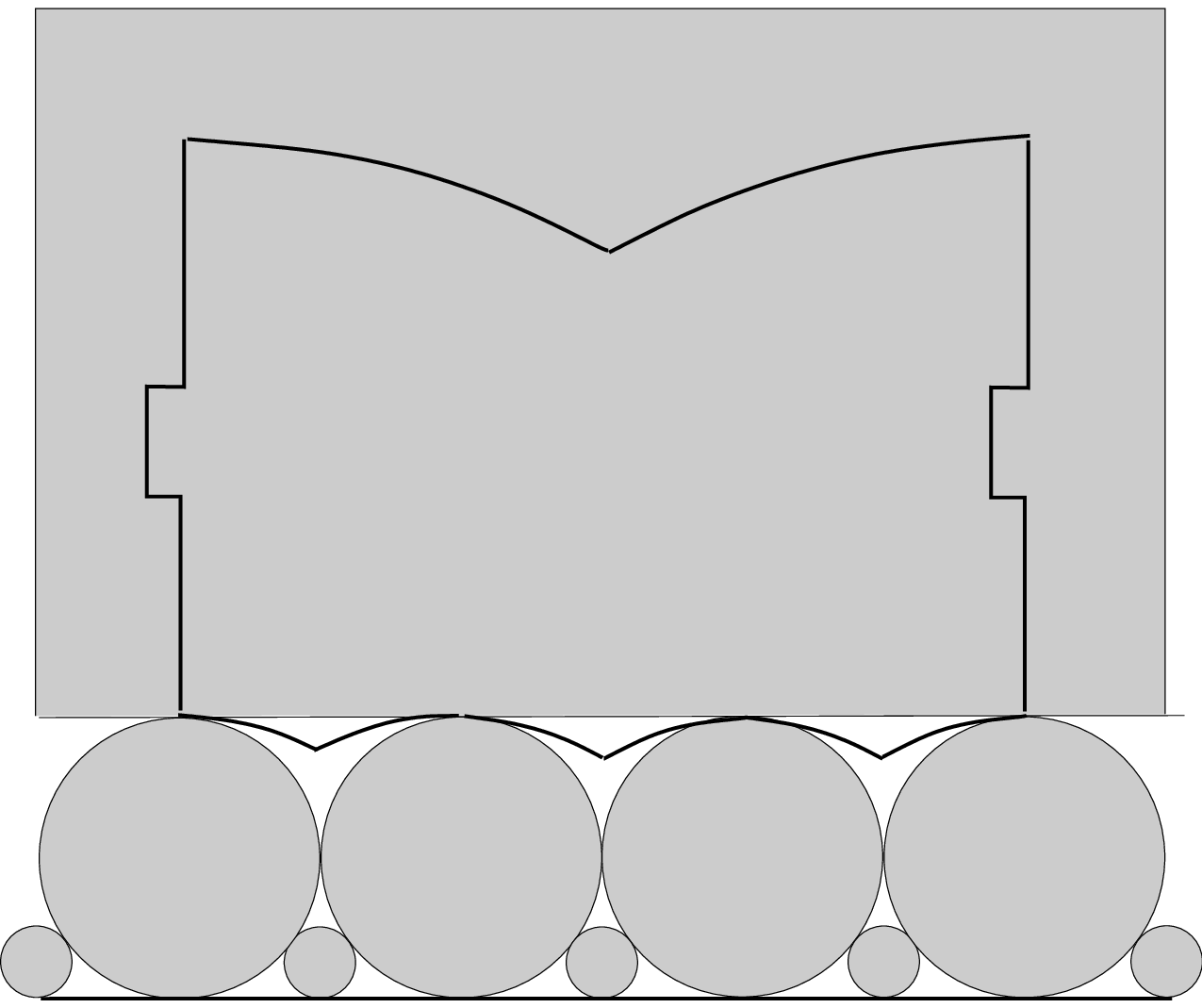;
\vs.3 \nd
\centerline{Figure 11. A triangular tile with 3 prongs}
\vfill\eject
\end